\documentclass[12pt]{article}

\usepackage{amsfonts}

\usepackage{amsmath,amssymb}

 \numberwithin{equation}{subsection}

\begin{document}

\title{Leading coefficients of the Kazhdan--Lusztig polynomials for
  an Affine Weyl group of type $\tilde{B}_2$}
\author{Liping Wang\\
\small Academy of Mathematics and Systems Science, \\
\small Chinese Academy of Sciences, Beijing, China;\\
\small China Economics and Management Academy,\\
\small Central University of Finance and Economics, Beijing, China\\
\small wanglp@amss.ac.cn}

\date{}

\maketitle

\begin{abstract}

In this paper we compute the leading coefficients $\mu (u,w)$ of the
Kazhdan--Lusztig polynomials $P_{u,w}$ for an affine Weyl group of
type $\tilde{B}_2$. By using the \textbf{a}-function of a Coxeter
group defined by Lusztig (see [L1, \S2]), we compute most $\mu
(u,w)$ explicitly. With part of these values $\mu (u,w)$, we show
that a conjecture of Lusztig on distinguished involutions is true
for an affine Weyl group of type $\tilde{B}_2$. We also show that
the conjectural formula in [L3, (12)] needs a modification.
\end{abstract}

\maketitle

\section*{}
The Kazhdan--Lusztig polynomials of a Coxeter group $W$ play a
 central role in Kazhdan--Lusztig theory. From [KL1, (2.2.c)] one sees that the
 leading coefficients $\mu (u,w)$ of some Kazhdan--Lusztig
 polynomials $P_{u,w}$ are  very important in understanding the Kazhdan--Lusztig
 polynomials. Moreover, the coefficients are of great importance in
 representation theory and Lie theory, and are related to some cohomology groups and
some difficult irreducible characters (see [A,CPS,S]).

In [L3] Lusztig computes the leading coefficients for some
Kazhdan--Lusztig polynomials of an affine Weyl group of type
$\tilde{B}_2$. In [S] for an affine Weyl group of type
$\tilde{A}_5$, some non-trivial leading coefficients are worked out.
McLarnan and Warrington  showed that $\mu (u,w)$ can be greater than
1 for a symmetric group (see [MW]). In [X1], Xi showed that when
$u\leq w$ and $a(u)< a(w),$ then $\mu (u,w)\leq 1$ if $W$ is a
symmetric group or an affine Weyl group of type $\tilde{A}_n.$ In
[SX, Theorem 3.3], Scott and Xi showed that the leading coefficient
$\mu (u,w)$ of some Kazhdan--Lusztig polynomial $P_{u,w}$ of an
affine Weyl group of type $\tilde{A}_n$ is $ n+2$ if $n\geq 4.$ In
[G], Green showed that when $W$ is a Coxeter group of type
$\tilde{A}_{n-1}$ ($n\geq3$), then $\mu(u,w)\leq1$ if $u$ is fully
commutative. In [J], Jones showed that $\mu(u,w)\leq1$ when $w$ is a
Deodhar element of a finite Weyl group.

In this paper we compute the coefficients $\mu (u,w)$ for an affine
Weyl group $W$ of type $\tilde{B}_2$.  There are four two-sided
cells in $W$: $c_e =\{w\in W\mid a(w)=0\}=\{e\},\ c_1=\{w\in W\mid
a(w)=1\},\ c_2=\{w\in W\mid a(w)=2\},\ \textrm{and}\ c_0=\{w\in
W\mid a(w)=4\}$, where $e$ is the neutral element of $W$ (see [L1,
\S11.2]). The main results are the following:

 (1) For any $u\le w$ in $W$ with
$(u,w)\not\in (c_0\times c_0)\cup (c_0\times c_1)\cup (c_0\times
c_2)$, the  value $\mu (u,w)$ is determined . The values of $\mu
(u,w)$ are displayed in Sections 3,  5,  7 and 8. When $(u,w)\in
c_0\times c_0$, the value $\mu (u,w)$ can  be computed by using a
formula in the proof for [SX,  Theorem 3.1], see [W] for details.

(2) Using the values of $\mu (u,w)$, we prove that Lusztig's
conjecture on distinguished involutions proposed in 1987
 is true for $W$ (see
Theorem 2.2) and that the  W-graph of type $\tilde{B}_2$ is
non-locally finite (see Theorem 8.2 and Remark 8.3).

(3) The conjectural formula for $b_{\lambda,\lambda'}$ in [L3, (12)]
needs a modification (see Proposition 5.1 and Remark 5.3.)

\def\L{\Lambda}
\def\l{\lambda}
\section{Preliminaries}In this section we recall some basic facts about $\mu (u,w)$ which will be needed later.

\subsection{Basic definitions and conventions}Let $R$ be a root system and $W_0$  its Weyl group.
Denote by $\Lambda$ the weight lattice of $R$ and $\Lambda_r=\mathbb
ZR$ the root lattice. The semi-direct product $W=W_0 \ltimes
\Lambda_r$ is an affine Weyl group and $\tilde{W}=W_0 \ltimes
\Lambda$ is  an extended affine Weyl group.

We shall denote by $S$ the set of simple reflections of $W.$ There
is  an abelian subgroup $\Omega$ of $\tilde{W}$ such that $\omega
S=S\omega$ for any $\omega \in\Omega$ and $\tilde{W}=\Omega \ltimes
W.$ The length function $l$ of $W$ and
 the Bruhat order $\le $on $W$  can be extended to $\tilde{W}$
by setting $l(\omega w)=l(w)$ and $\omega w\leq \omega' u$ if and
only if $\omega=\omega'$ and $w\leq u,$ where $\omega, \omega'$ are
in $\Omega$ and $w,u$ are in $W.$

\def\th{\tilde{\cal H}}
Let $\tilde{\cal H}$ be the Hecke algebra of $(\tilde W,S)$ over
$\mathcal {A}=\mathbb{Z}$ $[q^{\frac{1}{2}}, q^{-\frac{1}{2}}]$ ($q$
an indeterminate) with parameter $q.$ By definition, $\th$ is a free
$\mathcal A$-module and has a basis $\{T_w\}_{w\in \tilde W}$, its
multiplication is defined by the relations $(T_s-q)(T_s+1)=0$ if
$s\in S$ and $T_wT_u=T_{wu}$ if $l(wu)=l(w)+l(u)$. Let $C_w
=q^{-\frac{l(w)}{2}}\sum_{u\le w} P_{u,w}T_u, w\in \tilde W$ be its
Kazhdan--Lusztig basis, where $P_{u,w}\in \mathbb{Z} [q]$ are the
Kazhdan--Lusztig polynomials. (Here the $C_w$ is just the $C_w'$ in
[KL1].) The degree of $P_{u,w}$ is less than or equal to
$\frac{1}{2} (l(w)-l(u)-1)$ if $u<w$ and $P_{w,w}=1.$ See [X3, 1.6]
for more details. The subalgebra $\cal H$ of $\th$ generated by all
$T_s\ (s\in S$) is the Hecke algebra of $(W,S)$. For $\omega$ in
$\Omega$ and $w$ in $W$, we have $C_{\omega w}=T_\omega C_w$. Thus
we have $P_{\omega u,\omega w}=P_{u,w}$ and $P_{\omega' u,\omega
w}=0$ for any different $\omega,\ \omega'$ in $\Omega$ and $u,w$ in
$W$. It is known that the coefficients of these polynomials are all
non-negative (see [KL2, Corollary 5.6]).

Write $P_{u,w}=\mu (u,w)q^{\frac{1}{2} (l(w)-l(u)-1)}+$ lower degree
terms. The coefficient $\mu (u,w)$ is very interesting, which can be
seen even from the recursive formula for Kazhdan--Lusztig
polynomials [KL1, (2.2.c)]. We call $\mu (u,w)$ the Kazhdan--Lusztig
 coefficient of $P_{u,w}.$ Write $u\prec w$ if
$u\le w$ and $\mu (u,w)\neq0.$ Define $\tilde\mu (u,w)=\mu (u,w)$ if
$u\le w$ and $\tilde\mu (u,w)=\mu (w,u)$ if $w\leq u.$

The following  properties for  Kazhdan--Lusztig polynomials are known (see [KL1]):\\
(a) $P_{u,w}=P_{u^{-1},w^{-1}}.$ In particular,
$\mu(u,w)=\mu(u^{-1},w^{-1})$.\\
(b) Let $u,w\in W,\ s\in S$ be such that $u<w,\ su>u,\ sw<w$. Then
$u\prec w$ if and only if $w=su.$ Moreover,
$\mu(u,w)=1$ in this case.\\
(c) Let $u,w\in W,\ s\in S$ be such that $u<w,\ us>u,\ ws<w$. Then
$u\prec w$ if and only if $w=us.$ Moreover, $\mu(u,w)=1$ in this
case.

We  refer to [KL1] for the definition of the preorders $\leq_L,
\leq_R, \leq_{LR}$ and of the equivalence relations $\sim_L, \sim_R,
\sim_{LR}$ on $W.$ The corresponding equivalence classes are called
left cells, right cells, two-sided cells of $W,$ respectively. For
$w\in W$, set $\mathcal {L}(w)=\{s\in S\mid sw\leq w\}$ and
$\mathcal {R}(w)=\{s\in S\mid
ws\leq w\}.$ Then we have (see [KL1, Prop. 2.4])\\
(d) $\mathcal {R}(w)\subseteq \mathcal {R}(u) \;\textrm{if}\;u\leq_L
w.$ In particular,
$\mathcal {R}(w)=\mathcal {R}(u)\;\textrm{if}\;u\sim_L w;$\\
(e) $\mathcal {L}(w)\subseteq \mathcal {L}(u)\;\textrm{if}\;u\leq_R
w.$ In particular, $\mathcal {L}(w)=\mathcal
{L}(u)\;\textrm{if}\;u\sim_R w$.

The $\textbf{a}$-function on $W$ is defined in [L1, \S2], which is a
useful tool to study cells of $W$. The following result is due to
Springer (see [X1, 1.3(a)])\\
(f) If $\tilde\mu (u,w)$ is non-zero, then (i) $a(u)<a(w)$ implies
$w\leq_L u$ and $w\leq_R u$; (ii) $a(u)=a(w)$ implies that $ u\sim_L
w$ or $u\sim_R w$.

Write $C_w C_u=\sum_{z\in W}h_{w,u,z}C_z,$ where $h_{w,u,z}\in
{\mathcal {A}}.$ Following Lusztig and Springer, we define
$\delta_{w,u,z}$ and $\gamma_{w,u,z}$ by the following formula:
$h_{w,u,z}=\gamma_{w,u,z} q^\frac{a(z)}{2} + \delta_{w,u,z}
q^\frac{a(z)-1}{2} +$ lower degree terms.

We shall need a result of Lusztig. Fix a subset $S'\subset S$
consisting of two elements $s_1,s_2$ such that $s_1s_2$ has order
$4$ and we denote by $W'$ the subgroup generated by $s_1,s_2$. Each
coset $W'w$ can be decomposed into four parts: one consists of the
unique element $u$ of minimal length, one consists of the unique
element $u'$ of maximal length, one consists of the three elements
$s_1u,s_2s_1u,s_1s_2s_1u$, and one consists of the three elements
$s_2u,s_1s_2u,s_2s_1s_2u$. The last two subsets are called (left)
$strings$. We shall regard them as sequences (as above) rather than
subsets. The following result is due to Lusztig (see [L1,
(10.4.2)]).\\
(g) Consider two strings $u_1,u_2,u_3$ and $w_1,w_2,w_{3}$ (with
respect to $S'$). Set $a_{ij}=\tilde{\mu}(u_i,w_j)$ if $S'\cap
\mathcal {L}(u_i)=S'\cap \mathcal {L}(w_j)$ and $a_{ij}=0$
otherwise. Then the integers $a_{ij}$ satisfy the identities:
$a_{11}=a_{33},a_{13}=a_{31},a_{22}=a_{11}+a_{13}$, and
$a_{12}=a_{21}=a_{23}=a_{32}.$

We can also define (right) strings and then have similar identities.


\subsection {The lowest two-sided cell}
In this subsection we collect some facts about the lowest two-sided
cell of the extended affine Weyl group $\tilde{W}$.

For any $u=\omega_1u_1,\ w=\omega_2w_1,\ \omega_1,\omega_2\in
\Omega,\ u_1,w_1\in {W},$ we  say that $u\leq_Lw$ (respectively
$u\leq_Rw$, respectively $u\leq_{LR}w$) if $u_1\leq_Lw_1$
(respectively
$\omega_1u_1\omega_1^{-1}\leq_R\omega_2w_1\omega_2^{-1}$,
respectively $u_1\leq_{LR}w_1$) (see [X3, \S1.11]).
 The left (respectively right, respectively two-sided) cells of
$\tilde{W}$ are defined as those of $W$. We also define $a(\omega
w)=a(w)$ for $\omega\in\Omega,
 w\in W$.

 It is known that (see [Sh1])
$c_0=\{w\in\tilde{W}\mid a(w)=l(w_0)\}$ is a two-sided cell, which
is the lowest one for the partial order $\leq_{LR},$ where $w_0$ is
the longest element of $W_0.$ We call $c_0$ the lowest two-sided
cell of $\tilde{W}.$

In [SX, \S2.1], the authors gave a description of $c_0$. Let $R^+$
(respectively $R^-$, respectively $\Delta$) be the set of positive
(respectively
 negative, respectively simple) roots
in the root system $R$ of $W_0.$ The set of dominant weights
$\Lambda^+$ is the set $\{z\in\Lambda\mid l(zw_0)=l(z)+l(w_0)\}.$
For each simple root $\alpha$ we denote by $s_\alpha$ the
corresponding simple reflection in $W_0$ and $x_\alpha$ the
corresponding fundamental weight. For each $w\in W_0,$ we set
$$d_w=w\prod_{{\alpha\in\Delta}\atop {w(\alpha)\in R^-}}x_\alpha.$$
Then $c_0=\{d_wzw_0d_u^{-1}\mid w,u\in W_0, z\in\Lambda^+\}.$

\noindent (a) An element $w$ is in $c_0$ if and only if $w=uw_0u'$
for some $u,u'\in \tilde W$ such that $l(uw_0u')=l(u)+l(w_0)+l(u')$,
see [Sh1].

\section{  A conjecture of Lusztig}

Let $\delta (z)$ be the degree of $P_{e,z},$ where $e$ is the
neutral element of $W.$ Then actually one has $l(z)-a(z)-2\delta
(z)\geq 0$ (see [L2]). Set ${\cal D}_i=\{d\in W\mid
l(d)-a(d)-2\delta (d)= i\}.$ The elements in ${\cal D}_0$ are called
distinguished involutions of $W$. In [L2], Lusztig showed that
${\cal D}_0$ is finite for affine Weyl groups.

Lusztig has a conjecture to describe ${\cal D}_1$ (given in 1987,
see [Sh2]).

\noindent {\bf Conjecture 2.1} Assume that $W$ is a Weyl group or an
affine Weyl group. Let $z\in W.$ Then $z$ is in ${\cal D}_1$ if and
only if there exists some $d\in {\cal D}_0$ such that $z\sim_{LR}d$
and $\tilde{\mu}(z,d)\neq0.$

The "only if" part of the conjecture is true, this is due to Shi;
see [Sh2]. Moreover, for $z\in W$ and $d\in{\cal D}_0,$ if
$\tilde{\mu}(z,d)\neq0,z\sim_{LR}d$ and $z\nsim_Lz^{-1},$ it is not
difficult to prove that $z\in{\cal D}_1.$ This part is due to
Springer; see [Sh2]. For a Weyl group $W_0$, Lusztig  showed that
Conjecture 2.1 holds whenever the complex representation of $W_0$
provided by the two-sided cell of $W_0$ containing $z$ does not
contain the irreducible representation of $W_0$ of degree 512 for
type $E_7$ or of degree 4096 for type $E_8$ (see [X1, Sh2]). Also Xi
proved that this conjecture is true for an affine Weyl group of type
$\tilde{A}_n$; see [X2, Theorem 1.5].

In section 6 we will prove the following result.\\
 {\bf Theorem 2.2 } The Conjecture 2.1 is
true for an affine Weyl group of type $\tilde{B}_2$.

\section {The second highest two-sided cell}
In this section, we assume that $(W,S)$ is an arbitrary affine Weyl
group. The second highest two-sided cell $c_1$ of $W$ is described
in [L4].

\noindent {\bf Proposition 3.1} We have $c_1=\{w\in W\mid
a(w)=1\}\\=\{e\neq w\in W\mid w\ \textrm{has a unique reduced
expression}\}.$

Then we get following corollary.

\noindent {\bf Corollary 3.2} \\
(i) For any $w\in c_1,\ |\mathcal {L}(w)|=|\mathcal {R}(w)|=1.$\\
(ii) Let $w\in c_1$. Assume that $\mathcal {L}(w)=\{s'\}$ and
$\mathcal {R}(w)=\{s''\},$
then $w\sim_Rs'$ and $w\sim_Ls''.$\\
(iii) For any $u,w$ in $c_1,$   $u\sim_Lw$ is equivalent to $
\mathcal {R}(u)=\mathcal {R}(w)$ and $u\sim_Rw$ is equivalent $
\mathcal {L}(u)=\mathcal {L}(w).$

 \noindent {\bf Proof.} We only need to
prove (ii). By  assumption
 we have  $w=w's''$ for some $w'\in
W$ and $l(w)=l(w')+1$.  Then we get $w\leq_Ls''.$ By [L2, Corollary
1.9(b)] and $a(w)=a(s'')=1,$ we get $w\sim_Ls''.$ Similarly, we get
that $w^{-1}\sim_Ls'$,  then $w\sim_Rs'$ follows.\hfill$\Box$\

Now we give the values of $\mu(u,w)$ for $u,w\in c_1$.

\noindent {\bf Proposition 3.3} Assume that $(W,S)$ is an
irreducible affine Weyl group of type different from
$\tilde{A}_{2}$. Then for any $u,w\in c_1,\ u\leq w, $ we have
$${\mu (u,w)}=\left\{\begin{array}{ll}
{1}, & \ \ \textrm{if} \;  l(w)-l(u)=1 ,\\
{0} ,   & \ \ \textrm{otherwise}.\end{array} \right.$$

\noindent {\bf Proof.} Let $u\leq w$ be in $c_1$. Assume that $\ \mu
(u,w)\neq 0.\ $ By 1.1 (f) we have  $\ u\sim_L w\ $ or $\ u\sim_R w\
$. If $\ u\sim_L w\ $ but  $\ u\nsim _R w, $ then by Corollary 3.2,
there exists some $s'\in S$ such that $s'u>u$ and $s'w<w$. By
1.1(b), we get that $w=s'u$ and  $\mu (u,w)=1$. Similarly we can get
that $l(w)-l(u)=1$ and $\mu(u,w)=1$ if $u\nsim_L w$ but $u\sim_R w$.

To complete the proof we need to prove that $\mu(u,w)=0$ for $u\le
w$ in $c_1$ satisfying $\mathcal {L}(u)=\mathcal {L}(w)\
\textrm{and}\ \mathcal {R}(u)=\mathcal {R}(w)$.

First we assume that $W$ is not of type $\tilde{A}_{2m}$, so that
the Coxeter graph contains no odd cycles. By Proposition 3.1, it can
be checked that $l(w)-l(u)$ is even for any $u,w$ in $ c_1$
satisfying $\mathcal {L}(u)=\mathcal {L}(w)\ \textrm{and}\ \mathcal
{R}(u)=\mathcal {R}(w)$. We are done in this case.

 When $W$ is of type $\tilde{A}_{2m}$
for $m\geq 2,$ the result can be proved by using star operation
defined in [KL1]. Let $S=\{s_0,s_1,\ldots,s_n\}$ be the Coxeter
generators, where $n=2m$. We know that they satisfy relations
$(s_is_{i+1})^3=1,$ for $0\leq i\leq n-1$, $(s_ns_0)^3=1$ and
$s_{i}s_{j}=s_js_i$ for $1\leq|j-i|\leq n-1$. Using symmetry and
star operations we are reduced to the case $u=(s_1s_2s_3\cdots
s_ns_0)^{l_1}$ for some $l_1\in \mathbb{N}$ and $w=(s_1s_2s_3\cdots
s_ns_0)^{k}$ or $w=(s_1s_0s_ns_{n-1}\cdots s_3s_2)^{k}s_1s_0$ for
some $k\in \mathbb{N}$ such that $l_1\leq k$, where $\mathbb N$ is
the set of nonnegative inters. We can assume that $k-l_1=2k+1$ for
some $k \in \mathbb{N}$. If $w=(s_1s_2s_3\cdots s_ns_0)^{k}$, then
by [KL1, Theorem 5.2], we get that
\begin{eqnarray*}
\mu(u,w)&=&\mu((s_1s_2s_3\cdots s_ns_0)^{l_1},(s_1s_2s_3\cdots s_ns_0)^{k}) \\
&=&\mu(s_1s_2s_3\cdots s_ns_0,(s_1s_2s_3\cdots s_ns_0)^{k-l_1+1})\\
&\vdots& \\
&=&\mu(s_0,s_0(s_1s_2s_3\cdots s_ns_0)^{k-l_1})\\
&=&\mu(s_1s_0,(s_1s_2s_3\cdots s_ns_0)^{k-l_1})\\
&\vdots& \\
&=&\mu(u',w'),
\end{eqnarray*}
where $u'=s_ms_{m-1}\cdots s_2s_1s_0(s_n\cdots s_2s_1s_0)^k\ $ and\\
$w'=s_ms_{m+1}\cdots s_ns_0(s_1s_2\cdots s_ns_0)^k$. We have
$l(w')-l(u')=1$ since $n=2m$. Because $m\geq2$ and $w$ has a unique
reduced expression, we have $u'\nleq w'$. Thus we obtain that
$\mu(u,w)=\mu(u',w')=0$ in this case.  If $w=(s_1s_0s_n\cdots
s_3s_2)^{k}s_1s_0$, we can prove that $\mu(u,w)=0$ similarly.

The proof is completed.\hfill$\Box$\

\noindent {\bf Proposition 3.4} Let $(W,S)$ be an irreducible affine
Weyl group of type $\tilde{A}_{2}$. Then for any $u,w\in c_1,\ u\leq
w, $
we have \\
(i) If $\mathcal {L}(u)\neq \mathcal {L}(w)$ or $\mathcal {R}(u)\neq
\mathcal {R}(w)$, then
$${\mu (u,w)}=\left\{\begin{array}{ll}
{1}, & \ \ \textrm{if} \;  l(w)-l(u)=1,\\
{0} ,   & \ \ \textrm{otherwise}.\end{array} \right.$$ (ii) If
$\mathcal {L}(u)= \mathcal {L}(w)$ and $\mathcal {R}(u)= \mathcal
{R}(w)$, then
$${\mu (u,w)}=\left\{\begin{array}{ll}
{1}, & \ \ \textrm{if} \;  l(w)-l(u)=3,\\
{0} ,   & \ \ \textrm{otherwise}.\end{array} \right.$$

 \noindent{\bf Proof.} Part (i) is obvious by Corollary 3.2 and 1.1 (b)--(c).
 Let $s_0,s_1,s_2$ be the simple reflections of $W$, so that $(s_0s_1)^3=(s_1s_2)^3=(s_0s_2)^3$.

Assume that $\mathcal {L}(u)=\mathcal {L}(w)\ \textrm{and}\ \mathcal
{R}(u)=\mathcal {R}(w)$. Without loss of generality, we can assume
that $u=(s_1s_2s_0)^mz$ for some $z\in\{e,s_1,s_1s_2\}$ and
$m\in\mathbb{N}$. Then $w=(s_1s_2s_0)^nz$ or $w=(s_1s_0s_2)^nv$ for
some $v\in\{e,s_1,s_1s_0\}$ and $n\in\mathbb{N}$. We  may require
that $n-m=2k+1$ for some nonnegative integer $k$ to ensure that
$l(w)-l(u)$ is odd and $u\le w$.

If $w=(s_1s_2s_0)^nz$, then by [KL1, Theorem 4.2], we have
\begin{eqnarray*}
\mu(u,w)&=&\mu((s_1s_2s_0)^m,(s_1s_2s_0)^n) \\
&=&\mu(s_1s_2s_0,(s_1s_2s_0)^{n-m+1})\\
&=& \mu(s_2s_0,s_2s_0(s_1s_2s_0)^{n-m})\\
&=&\mu(s_0,s_0(s_1s_2s_0)^{n-m})\\
&=&\mu(s_1s_0,(s_1s_2s_0)^{n-m})\\
&\vdots& \\
&=&\mu(s_0(s_2s_1s_0)^i,s_0(s_1s_2s_0)^{n-m-i})\\
&=&\mu(s_1s_0(s_2s_1s_0)^i,(s_1s_2s_0)^{n-m-i})\\
&=&\mu((s_2s_1s_0)^{i+1},s_2s_0(s_1s_2s_0)^{n-m-i-1})\\
&\vdots& \\
&=&\mu(s_1s_0(s_2s_1s_0)^k,(s_1s_2s_0)^{k+1}).
\end{eqnarray*}
If $n-m=1$, then $k=0$ and $\mu(u,w)=\mu(s_1s_0,s_1s_2s_0)=1$. If
$n-m>1$, then $k>0$ and $s_1s_0(s_2s_1s_0)^k\nleq
(s_1s_2s_0)^{k+1}$. So in this case we have
$\mu(u,w)=\mu(s_1s_0(s_2s_1s_0)^k,(s_1s_2s_0)^{k+1})=0$.

Now suppose that $w=(s_1s_0s_2)^nv$ for some $v\in\{e,s_1,s_1s_0\}$.
We can assume that $z=e$, so that $v=s_1s_0$. If $k>0$, then by
[KL1, Theorem 4.2], we have
\begin{eqnarray*}
\mu(u,w)&=&\tilde{\mu}((s_1s_2s_0)^m,(s_1s_0s_2)^ns_1s_0) \\
&=&\tilde{\mu}(s_0(s_1s_2s_0)^m,s_0s_2(s_1s_0s_2)^{n-1}s_1s_0)\\
&=& \tilde{\mu}(s_2s_0(s_1s_2s_0)^m,s_2(s_1s_0s_2)^{n-1}s_1s_0)\\
&=&\tilde{\mu}((s_1s_2s_0)^{m+1},(s_1s_0s_2)^{n-1}s_1s_0).
\end{eqnarray*}
This reduces the problem to the case $k=0$, and we can prove that
$\mu(u,w)=0$ in this case because $u\nleq w$.

Thus (ii) is true. The proof is completed.
 \hfill$\Box$\
\section {Semilinear equations related to $\mu$}

 In [L3] Lusztig introduced some semilinear equations, which are useful for calculating some $\mu(u,w)$.
 In this section we first recall the equations, then give some discussions to the type $\tilde{B}_2$.

\def\L{\Lambda}
\def\l{\lambda}

 We need some notations. Let $R^+$ be the subset of $R$ containing all positive roots. For any subset
$\textbf{i}\subseteq R^+,$ we set
$\alpha_{\textbf{i}}=\sum_{\alpha\in \textbf{i}}\alpha \in \Lambda$.
Set $\rho=\frac 12\alpha_{R^+}.$ Then $\rho$ is a dominant weight
and $w(\rho)-\rho\in \Lambda_r$ for any $w\in W_0.$ For any
$\lambda\in \Lambda_r$ we set
$$\Phi(\lambda)=\sum_{{}\atop{\textbf{i}\subseteq R^+};\
\alpha_{\textbf{i}} =\lambda}(-v^2)^{-|\textbf{i}|}.$$ Note that the
summation index $\textbf{i}$ runs through all subsets $\textbf{i}$
of $R^+$ such that $\lambda$ can be written as the sum of all
elements in $\textbf{i}$. If $0\neq\lambda$ cannot be written as a
sum of distinct positive roots,  we set $\Phi(\lambda)=0.$

Let $\L_r^+$ be the set of dominate weights in the root lattice
$\L_r$, that is $\Lambda_r^+=\Lambda_r\cap \Lambda^+.$ For any
$\lambda\in \Lambda_r^+,$ we set $$W_0^{\lambda}=\{w\in W_0\mid
w(\lambda)=\lambda\}$$ and
$$\pi_{\lambda}=v^{-\nu_{\lambda}}\sum_{w\in
W_0^{\lambda}}v^{2l(w)},$$ where $\nu_{\lambda}$ is the number of
reflections of $W_0^{\lambda}.$ For  $\lambda$ and $\lambda'$  in $
\Lambda_r^+$ we define
$$a_{\lambda,\lambda'}=\frac{v^{\nu_{\lambda'}}}{\pi_{\lambda'}}\sum_{w\in W_0}
(-1)^{l(w)}\Phi(\lambda'+\rho-w(\lambda+\rho)).$$

 \noindent {\bf Convention:} For any element in $\Lambda$, we will
use the same notation when it is regarded as an element in
$\tilde{W}$. For two elements in $\L$, the operation between them
will be written additively if they are regarded as elements in
$\Lambda$ and will be written  multiplicatively when they are
regarded as elements in $\tilde{W}$.

Lusztig showed that there is a 1-1 correspondence between
$\Lambda_r^+$ and the set of $W_0-W_0$ double cosets in $W$ (an
element $\lambda$ of $\Lambda_r^+$ corresponds to the unique double
coset $W_0\lambda W_0$ containing it: see [L5, \S2]). For each
$\lambda\in \Lambda_r^+,$ there is a unique element $m_\lambda$ of
minimal length and a unique element $M_\lambda$ of maximal length in
$W_0\lambda W_0.$ We have $\lambda\leq\lambda'$ (i.e. $\l'-\l\in
\mathbb NR^+$) if and only if $M_\lambda\leq M_{\lambda'}$ (Bruhat
order) for $\lambda,\lambda'\in \Lambda_r^+$.

Set $v=q^{\frac{1}{2}},\ p_{u,w}=v^{l(u)-l(w)}P_{u,w}(v^{2})\in
\mathbb{Z} [v^{-1}]$ for $u\leq w\in W$ and $p_{u,w}=0$ for all
other $u,w$ in $W.$ For two elements $\l,\ \l'$ in $\L_r^+$, set
$$b_{\lambda,\lambda'}=\sum_{z\in W_0\lambda
W_0}(-v)^{l(m_\lambda)-l(z)}p_{z,m_{\lambda'}}.$$ By [L3, \S5], we
know

\noindent (a)   $a_{\lambda,\lambda'}$ is zero unless
$\lambda\leq\lambda'$ and are equal to 1 when $\lambda=\lambda';$
when $\lambda<\lambda'$ it belongs to $v^{-1}\mathbb{Z}[v^{-1}].$

\noindent (b) $b_{\lambda,\lambda'}$ is zero unless
$\lambda\leq\lambda'$ and is equal to 1 when $\lambda=\lambda';$
when $\lambda<\lambda'$ it belongs to $v^{-1}\mathbb{Z} [v^{-1}].$
Also we have
$\mu(m_\lambda,m_{\lambda'})=\textrm{Res}_{v=0}(b_{\lambda,\lambda'}),$
where $\textrm{Res}_{v=0}(f)\in \mathbb{Z} $ denotes the coefficient
of $v^{-1}$ in $f\in {\mathcal {A}}=\mathbb Z[v,v^{-1}].$

Let $\bar{}:\ {\mathcal {A}}\rightarrow{\mathcal {A}} $ be the ring
involution such that $\bar{v}=v^{-1}.$ The following lemma of
Lusztig gives a way to compute $b_{\l,\l'}$ inductively (see [L3,
Proposition 7]).

\noindent {\bf Lemma 4.1} For any $\lambda,\lambda''\in
\Lambda_r^+,$ we have
$$\sum_{\lambda'\in
\Lambda_r^+}a_{\lambda,\lambda'}(-1)^{l(m_{\lambda'})-l(M_{\lambda'})}
\pi_{\lambda'}\bar{b}_{\lambda',\lambda''}=\sum_{\lambda'\in
\Lambda_r^+}\bar{a}_{\lambda,\lambda'}(-1)^{l(m_{\lambda'})-l(M_{\lambda'})}
\pi_{\lambda'}b_{\lambda',\lambda''}.$$

\def\L{\Lambda}
\def\l{\lambda}
\def\a{\alpha}
\def\b{\beta}
\def\a{\alpha}
\def\b{\beta}

In the rest of this paper,  $R$ is a root system of type $B_2$. Then
$(W,S)$ is an affine Weyl group of type $\tilde{B}_2$ and
$\tilde{W}$ is an extended affine Weyl group associated with $R$.

Let $\a$ and $\b$ be the long and short simple root of $R$
respectively. Then the set $R^+$ of positive roots in $R$ consists
of the following four elements: $ \alpha,$ $ \beta,\alpha+\beta,$ $
\alpha+2\beta .$ The fundamental dominant weights are
$$x=\alpha+\beta\quad{\text{ and }}\quad
y=\frac{1}{2}\alpha+\beta.$$ So the set $\L^+$ of dominant weights
consists of the elements $mx+ny,\ m,n\in\mathbb N.$ Also
 we have
$$\Lambda_r^{+}=\{mx+2ny \mid m,n\in
\mathbb{N}\}=\{i\alpha+j\beta\mid  i,j\in \mathbb{N},i\leq
j\leq2i\}.$$

Let $s$ and $t$ be the simple reflections in $S$ corresponding to
$\a$ and $\b$ respectively. The Weyl group $W_0$ is generated by $s$
and $t$. The set $S$ contains a unique element
 out of $W_0$, denoted by $r$. Then $W$ is generated by $r,s,t$. We have
  $rt=tr,\ (rs)^4
=(st)^4 =e.$

Define $X_1$ (respectively $X_2$)  to be the subsets of
$\Lambda_r^+$ consisting of elements $mx,\ m\geq1)$ (respectively
$2my,$ $m\ge 1$). Let $Y_1$ (respectively $Y_2$)  be the subsets of
$\Lambda_r^+$ consisting of elements $x+2my, \ m\geq1)$
(respectively $nx+2my,$ $n\ge 2$ and $ m\ge 1$). Set $X=X_1\cup
X_2,\ Y=Y_1\cup Y_2$. We have $\Lambda_r^{+}=X\cup Y\cup\{0\}$.

The following result is needed in calculating $a_{\l,\l'}$ and
$b_{\l',\l}$.

\noindent {\bf Proposition 4.2} Let $\lambda$ be in $X\cup Y$. We have\\
(i) If $\l\in X_1$ then $W_{0}^{\lambda}= \{e,t\}.$ If $\l\in X_2$
then $W_0^\l=\{e,s\}$. If $\l\in Y$ then $W_0^\l=\{e\}.$\\
 (ii)
$\nu_{\lambda}=\left\{\begin{array}{ll}
1,& \ \ \textrm{if} \ \lambda\in X,\\
0,& \ \ \textrm{if} \ \lambda\in Y.\end{array} \right.$\\
(iii) $\pi_{\lambda}=\left\{\begin{array}{ll}
{v+v^{-1}},& \ \ \textrm{if} \ \lambda\in X,\\
1,& \ \ \textrm{if} \ \lambda\in Y.\end{array} \right.$\\
(iv) $${m_{\lambda}}=\left\{\begin{array}{ll}
r(stsr)^{m-1},& \ \ \textrm{if} \ \lambda=mx\in X_1,\\
rsr(tsr)^{2m-2},& \ \ \textrm{if} \ \lambda=2my\in X_2,\\
rsr(tsr)^{2m-1},& \ \ \textrm{if} \ \lambda=x+2my\in Y_1,\\
(rst)^{2m-1}(srst)^{n-2}srsrtsr,& \ \ \textrm{if} \
\lambda=nx+2my\in Y_2.\end{array} \right.$$ (v)
$(-1)^{l(m_{\lambda})-l(M_{\lambda})}=\left\{\begin{array}{ll}
-1,& \ \ \textrm{if} \ \lambda\in X,\\
1,& \ \ \textrm{if} \ \lambda\in Y.\end{array} \right.$\\

\noindent {\bf Proof.} For the action of a Weyl group on the
corresponding weight lattice, a  well known result says that the
stabilizer of a dominant weight is generated by the simple
reflections stablizing the dominant weight. The simple reflections
in $W_0$ are $s,t$ and we have $s(x)=x-\alpha=2y-x,\ s(y)=y,\
t(x)=x,\ t(y)=y-\beta=x-y$. Thus (i) follows from these facts. Also
(i) can be proved by a direct calculation since $W_0$ only contains
8 elements.

\def\ww{\tilde{W}}
(ii) is just a consequence of (i), since $\nu_{\lambda}$ is the
number of reflections in $W_{0}^{\lambda}$. We know that
$\pi_{\lambda}=v^{-\nu_{\lambda}}\sum_{w\in
W_0^{\lambda}}v^{2l(w)}$, thus (iii) can be checked by (i) and (ii).

Now we prove (iv). When $\l\in\L_r^+$, in $\ww$ we have $\l
w_0=w_0\l^{-1}$ and $l(\l w_0)=l(w_0)+l(\l)$. For $\l\in\L_r^+$,
there is a unique $w\in W_0$ such that (1) $\l=wm_\l$, (2)
$l(\l)=l(w)+l(m_\l)$, (3) $\mathcal L(m_\l)=\{r\}$.

In $\tilde{W}$ we have $x=stsr$ and $y=\omega rsr$, where $\omega\in
\Omega$ satisfies
 $\omega t=r\omega$. Thus
$y^2=tstrsr=tsrtsr$.

If $\lambda=mx \in X_1$, then in $\tilde{W}$ we have
$\lambda=(stsr)^m=stsr(stsr)^{m-1}$. Therefore
$m_{\lambda}=r(stsr)^{m-1}$. If $\lambda=2my\in X_2$, we get that
$\lambda=(tstrsr)^{m}=tstrsr(tstrsr)^{m-1}$. So
$m_{\lambda}=rsr(tstrsr)^{m-1}=rsr(tsr)^{2m-2}$, since $tr=rt$.

\def\ww{\tilde{W}}

The longest element in $W_0$ is $w_0=stst$. When $\lambda=x+2my\in
Y_1$, in $\ww$ we have $\l=stsr(tstrsr)^m=stsrtstrsr(tstrsr)^{m-1}$.
Since $tr=rt$, we get
$\l=w_0rsrtsr(tsrtsr)^{m-1}=w_0rsr(tsr)^{2m-1}$. Hence
$m_\l=rsr(tsr)^{2m-1}$.

Now assume that $\lambda=nx+2my\in Y_2$. Set $\l_1=x+2my$,
$\l_2=(n-1)x$. Then in $\ww$ we have
$\l=\l_1\l_2=w_0rsr(tsr)^{2m-1}(stsr)^{n-1}$. Thus
$m_\l=rsr(tsr)^{2m-1}(stsr)^{n-1}$. Noting that $rsrs=srsr$,
$tsts=stst$ and $tr=rt$, we see
$m_{\lambda}=(rst)^{2m-1}(srst)^{n-2}srsrtsr$. Thus (iv) holds.

 Since for  $\l\in\L_r^+$, in $\ww$ we have $\l
w_0=w_0\l^{-1}$ and $l(\l w_0)=l(w_0)+l(\l)$, so
$M_{\lambda}=\lambda stst$. Thus  we have
$M_{\lambda}=stsm_{\lambda}stst$ if $\lambda\in X_1$,
$M_{\lambda}=tstm_{\lambda}stst$ if $\lambda\in X_2$. Also
$M_{\lambda}=ststm_{\lambda}stst$ if $\lambda\in Y$ since
$\l=w_0m_\l$ in this case. Hence (v) holds. The proof is complete.
\hfill$\Box$\

A simple computation leads to the following identities:\\
$\Phi(0)=1,$\\
$\Phi(\alpha)=\Phi(\beta)=-v^{-2},$\\
$\Phi(\alpha+\beta)=\Phi(\alpha+2\beta)=v^{-4}-v^{-2},$\\
$\Phi(2\alpha+2\beta)=\Phi(2\alpha+3\beta)=v^{-4}-v^{-6},$\\
$\Phi(2\alpha+\beta)=\Phi(\alpha+3\beta)=v^{-4},$\\
$\Phi(3\alpha+3\beta)=\Phi(2\alpha+4\beta)=-v^{-6},$\\
$\Phi(3\alpha+4\beta)=v^{-8}$ and $\Phi(\lambda)=0$ for all for all
other $\lambda$ in $\Lambda_r$.

Using the above formulas we can compute $a_{\l,\l'}$, which are
needed in determining $b_{\l,\l'}$ next section.

\noindent {\bf Proposition 4.3} Let $\lambda$ and $\lambda'$ be in
$\Lambda_r^{+}$. Assume that $0<\l<\l'$ and $\lambda'=i\alpha+j\beta$ ( $i\leq j\leq2i$). \\
(i) If
$\lambda'-\lambda\neq\alpha,\beta,\alpha+\beta,\alpha+2\beta,\alpha+3\beta,
2\alpha+\beta,2\alpha+2\beta,2\alpha+3\beta,2\alpha+4\beta,
3\alpha+3\beta$ or $3\alpha+4\beta$, we have
$a_{\lambda,\lambda'}=0$.\\
(ii) If $\lambda'-\lambda=\alpha$ or $\beta$, we have $a_{\lambda,\lambda'}=-v^{-2}$.\\
(iii) If $\lambda'-\lambda=\alpha+\beta$, we have
$a_{\lambda,\lambda'}=\left\{\begin{array}{ll}
0,& \ \ \textrm{when} \ j=i,\\
-v^{-2},& \ \ \textrm{when} \ j=2i-1,\\
v^{-4}-v^{-2},& \ \ \textrm{otherwise}.\end{array} \right.$\\
(iv) If $\lambda'-\lambda=\alpha+2\beta$, we have
$a_{\lambda,\lambda'}=\left\{\begin{array}{ll}
-v^{-2},& \ \ \textrm{when} \ j=i+1,\\
v^{-4}-v^{-2},& \ \ \textrm{when}\ j>i+1.\end{array} \right.$\\
(v) If $\lambda'-\lambda=\alpha+3\beta$ or $2\alpha+\beta$, we have $a_{\lambda,\lambda'}=v^{-4}$.\\
(vi) If $\lambda'-\lambda=2\alpha+2\beta$, we have
$a_{\lambda,\lambda'}=\left\{\begin{array}{ll}
0,& \ \ \textrm{when} \ j>i,\\
v^{-4}-v^{-6},& \ \ \textrm{when}\ j=i.\end{array} \right.$\\
(vii) If $\lambda'-\lambda=2\alpha+3\beta$, we have\\
$$a_{\lambda,\lambda'}=\left\{\begin{array}{ll}
v^{-4},& \ \ \textrm{when} \ j=2i-1\ \textrm{or}\ j=i+1,\\
v^{-4}-v^{-6},& \ \ \textrm{otherwise}.\end{array} \right.$$\\
(viii) If $\lambda'-\lambda=2\alpha+4\beta$ or $3\alpha+3\beta$, we have $a_{\lambda,\lambda'}=-v^{-6}$.\\
(ix) If $\lambda'-\lambda=3\alpha+4\beta$, we have
$a_{\lambda,\lambda'}=v^{-8}$.

\noindent {\bf Proof.} By definition, $\rho$ is the half of the sum
of positive roots in $R^+$. So we have $\rho=\frac32\alpha+2\beta$.
Thus $\rho-s(\rho)=\alpha,$ $\rho-t(\rho)=\beta,$
$\rho-ts(\rho)=\alpha+3\beta,$ $ \rho-st(\rho)=2\alpha+\beta,$
$\rho-sts(\rho)=3\alpha+3\beta,$ $\rho-tst(\rho)=2\alpha+4\beta$ and
$\rho-stst(\rho)=3\alpha+4\beta$, since $s(\alpha)=-\alpha,$
$s(\beta)=\alpha+\beta,$ $t(\alpha)=\alpha+2\beta$ and
$t(\beta)=-\beta$. Since $0<\lambda\in\Lambda_r^{+}$, it is easy to
see that $sts(\lambda),tst(\lambda)$ and $stst(\lambda)=-\l$ are all
less than $0$. For a weight $\nu=a\a+b\b$ in $\L_r^+$, set
$h(\nu)=a+b$. Then we have
$$h(\lambda'+\rho-w(\lambda+\rho))=h(\l'-w(\l)+\rho-w(\rho))\ge 8$$ if  $w=sts,tst$ or $stst$. As a consequence, we have
$\Phi(\lambda'+\rho-w(\lambda+\rho))=0$ when $w=sts,tst$ or $stst$.

Now we are ready to compute $a_{\l,\l'}$.

(i) If $\lambda'-\lambda=n\beta$ for $n\geq2$, we have
$\lambda=i\alpha+(j-n)\beta$. Then  $i+n\leq j\leq2i$  since
$\l,\l'\in\L_r^+$ and $0<\l<\l'$. Clearly we have
$\Phi(\lambda'-\lambda)=0$. Note that
$$\lambda'+\rho-s(\lambda+\rho)=(2i-j+n+1)\alpha+n\beta,$$
$$\lambda'+\rho-t(\lambda+\rho)=(2j-2i-n+1)\beta,$$
$$\lambda'+\rho-ts(\lambda+\rho)=(2i-j+n+1)\alpha+(2i+n+3)\beta$$ and
$$\lambda'+\rho-st(\lambda+\rho)=(j-n+2)\alpha+(2j-2i-n+1)\beta.$$
We have $2i-j+n+1\geq 3$ and $2i-j+n+1>n$, so
$\Phi(\lambda'+\rho-s(\lambda+\rho))=0$. Since $2j-2i-n+1\geq
n+1\geq3,$ we get $ \Phi(\lambda'+\rho-t(\lambda+\rho))=0$. From
$2i+n+3\geq3n+3\geq9,$ one sees that
$\Phi(\lambda'+\rho-ts(\lambda+\rho))=0$. Finally $j-n+2\geq
n+2\geq4$ (note that $i\ge n$) implies that
$\Phi(\lambda'+\rho-st(\lambda+\rho))=0$.\

Therefore $a_{\lambda,\lambda'}=0$ if $\lambda'-\lambda=n\beta$ for
$n\geq2$. Similarly we have $a_{\l,\l'}=0$ if $\l'-\l=n\a$ for $n\ge
2$.

Also, by a similar computation we see $a_{\l,\l'}=0$ if $\l'-\l$ is
one of the following elements: $n\a+\b\ (n\ge 3)$, $\a+n\b\ (n\ge
4),$ $n\a+2\b\ (n\ge 3),$ $2\a+n\b\ (n\ge 5),$ $n\a+3\b \ (n\ge 4)$,
$3\a+n\b\ (n\ge 5)$, $n\a+4\b\ (n\ge 4)$, $m\a+n\b\ (m\ge 4,\ n\ge
5)$. Part (i) is proved.

(ii) If $\lambda'-\lambda=\beta$, then $\lambda=i\alpha+(j-1)\beta$.
Since $\l,\l'\in\L_r^+$ and $0<\l<\l'$, we have $2\le i+1\leq
j\leq2i$ and $i\geq1$. Thus
$$\lambda'+\rho-s(\lambda+\rho)=(2i-j+2)\alpha+\beta,$$
$$\lambda'+\rho-t(\lambda+\rho)=(2j-2i)\beta,$$  $$
\lambda'+\rho-ts(\lambda+\rho)=(2i-j+2)\alpha+(2i+4)\beta$$ and
$$\lambda'+\rho-st(\lambda+\rho)=(j+1)\alpha+(2j-2i)\beta.$$ So we get
$\Phi(\lambda'-\lambda)=-v^{-2},$
$\Phi(\lambda'+\rho-s(\lambda+\rho))=\Phi(2\alpha+\beta)=v^{-4}$ if
$j=2i$ and
$\Phi(\lambda'+\rho-s(\lambda+\rho))=\Phi((2i-j+2)\alpha+\beta)=0$
 if $j<2i$, and
$\Phi(\lambda'+\rho-t(\lambda+\rho))=\Phi(\lambda'+\rho-ts(\lambda+\rho))=\Phi(\lambda'+\rho-st(\lambda+\rho))=0$.

Thus  we get
$a_{\lambda,\lambda'}=\frac{v}{v+v^{-1}}(-v^{-2}-v^{-4})=-v^{-2}$
when  $j=2i$ and $a_{\lambda,\lambda'}=-v^{-2}$ when $j<2i$. That
is, $a_{\lambda,\lambda'}=-v^{-2}$ if $\lambda'-\lambda=\beta$.

Similarly we get $a_{\lambda,\lambda'}=-v^{-2}$ if
$\lambda'-\lambda=\a$.

Part (ii) is proved. For parts (iii)--(ix), the arguments are
similar. \hfill$\Box$\

In [L1, \S11.2], Lusztig described the left cells and two-sided
cells of $(W,S).$ For any subset $J$ of $S=\{r, s, t\}$, we denote
by $W^J$ the set of all $w\in W$ such that $\mathcal {R}(w)=J.$ In
the following we use the elements in $J$ to denote $J$. Then $(W,S)$
has $16$ left cells:
$A_{rs}=W^{rs},\;A_{rt}=A_{rs}t,\;A_s=A_{rt}s,\;A_r=A_sr,$
$A_{st}=W^{st},\;A_{rt}'=A_{st}r,\;A_s'=A_{rt}'s,\;A_t=A_s't,$
$B_{rt}=W^{rt}-(A_{rt}\cup A_{rt}'),\;B_s=B_{rt}s,\;B_r=B_sr,$
$B_t=B_st,\;C_r=W^r-(A_r\cup B_r),\;C_t=W^t-(A_t\cup B_t),$
$C_s=W^s-(A_s\cup A_s'\cup B_s),\;D_\emptyset=W^\emptyset=\{e\}.$

Set $c_e=D_\emptyset,\ c_1=C_r\cup C_s\cup C_t,\ c_2=B_r\cup B_s\cup
B_t\cup B_{rt}$, and $c_0=A_r\cup A_s\cup A_s'\cup A_t\cup
A_{rs}\cup A_{st}\cup A_{rt}\cup A_{rt}'.$

From [L1, \S11.2], we know that $c_e,\ c_1,\ c_2,\ c_0$ exhaust
two-sided cells of $W.$ We have $c_e =\{w\in W\mid a(w)=0\}=\{e\},\
c_1=\{w\in W\mid a(w)=1\},\ c_2=\{w\in W\mid a(w)=2\}$, and
$c_0=\{w\in W\mid a(w)=4\}.$

When $u,w$ in $c_1$ with $u\le w$, the value $\mu(u,w)$ is given in
Proposition 3.3. In the following sections we compute $\mu(u,w)$ for
other pairs $(u,w)$ except those in $(c_0\times c_0)\cup(c_0\times
c_1)\cup (c_0\times c_2)$.


\section {Computing $\mu(u,w)$ for $(u,w)\in c_2\times c_2$}

 Set $U=\{e,s,ts,rs\}$
and $V=\{e,s,st,sr\}$, then $c_2=\{urt(srt)^mv\mid u\in U,v\in
V,m\in \mathbb{N} \}.$ We will  compute $\mu(u,w)$ for $u,w\in c_2$
by means of the semilinear equations in Section 4. To do this we
first compute the $b_{\lambda,\lambda''}$ for
$\lambda,\lambda''\in\Lambda_r^+$. Our main results in this section
are Theorems 5.7,  5.8 and  5.9.

\noindent{(a)} Let $m$ be a positive integer and $\l_i=m\a+(m+i)\b$,
$i=0,1,...,m-1$. Then we have

$${b_{\lambda_i,\lambda_{m-1}}}=\left\{\begin{array}{ll}
{1},& \ \ \textrm{if} \ i=m-1,\\
{v^{-1}},& \ \ \textrm{if} \ i=m-2=0,\\
{-v^{-2}},& \ \ \textrm{if} \ i=m-2\geq1,\\
{0},& \ \ \textrm{if }0\le i\le m-3.\end{array} \right.$$

\noindent{\bf Proof.} By 4 (b) we  have $b_{\l_{m-1},\l_{m-1}}=1$.
 By Proposition 4.3 (ii-i),
we have

\noindent (1)  $a_{\l_i,\l_{i+1}}=-v^{-2}$ for $0\le i\le m-2$ and
$a_{\l_i,\l_j}=0$ if $j-i\ge 2$.

Since $\l_0=mx$ is in $X$ and $\l_i=(m-i)x+2iy$ is in $Y$ for
$i=1,2,...,m-1$, by Proposition 4.2 we see

\noindent (2) $(-1)^{l(m_{\l_i})-l(M_{\l_i})}$ is $-1$ if $i=0$ and
is $1$ if $i\ge 1$,

\noindent (3) $\pi_{\l_0}=v+v^{-1}$  and $\pi_{\l_i}=1$ if $i\ge 1$.

Now using Lemma 4.1 for $\l_i \ (0\le i\le m-2),\ \l_{m-1}$, and  4
(a), we get

\noindent (4) $\bar b_{\l_i,\l_{m-1}}-v^{-2}
\bar{b}_{\lambda_{i+1},\lambda_{m-1}}= b_{\l_i,\l_{m-1}}-v^{2}
{b}_{\lambda_{i+1},\lambda_{m-1}}$ if $1\le i\le m-2$ and

\noindent (5) $-(v+v^{-1})\bar b_{\l_0,\l_{m-1}}-v^{-2}
\bar{b}_{\lambda_1,\lambda_{m-1}}=-(v+v^{-1})
b_{\l_0,\l_{m-1}}-v^{2} {b}_{\lambda_1,\lambda_{m-1}}.$

Assume $m=2$. Since $b_{\l_1,\l_1}=1$, by (5)  we have
$-(v+v^{-1})\bar b_{\l_0,\l_{1}}-v^{-2} =-(v+v^{-1})
b_{\l_0,\l_{1}}-v^{2}.$ So
$\bar{b}_{\lambda_0,\lambda_1}-v=b_{\lambda_0,\lambda_1}-v^{-1}$,
which follows that $b_{\lambda_0,\lambda_1}=v^{-1}$ by 4.2 (b).

Now assume that $m\geq3$.  By (4) and 4 (b) we get
$b_{\l_{m-2},\l_{m-1}}=-v^{-2}$.  Applying 4 (b),  (4) when
$i=m-3\ge 1$ and (5) when $m=3$ we get that
$b_{\l_{m-3},\l_{m-1}}=0$. Now using induction on $i$, (4) and (5)
we see easily that $b_{\l_i,\l_{m-1}}=0$ if $0\le i\le m-3$. This is
the last part of (a). \hfill$\Box$

\def\g{\gamma}
\noindent (b)  Let $m\ge 2$ be a positive integer. Set
$\l_i=m\a+(m+i)\b$
 and $\gamma_i=(m-1)\a+(m-1+i)\b$
for $i=0,1,...,m-1$.  Then we have
$${b_{\g_i,\lambda_{m-1}}}=\left\{\begin{array}{ll}
{v^{-1}},& \ \ \textrm{if} \ i=m-1\ge 2\ \textrm{or}\ i=m-2=0,\\
{v^{-4}-v^{-2}},& \ \ \textrm{if} \ i=m-2\ge 1,\\
{-v^{-3}},& \ \ \textrm{if} \ i=m-3=0,\\
{v^{-4}},& \ \ \textrm{if} \ i=m-3\ge 1,\\
{0},& \ \ \textrm{if }  i=m-1=1\ \textrm{or if }0\le i\le
m-4.\end{array} \right.$$

\noindent{\bf Proof.}  By Proposition 4.3, we have

\noindent (1)
$a_{\gamma_i,\gamma_{i+1}}=a_{\gamma_{i+1},\l_i}=-v^{-2}$ for $0\le
i\le m-2$ and $a_{\gamma_i,\gamma_j}=0$ if $j-i\ge 2$;

\noindent (2) ${a_{\gamma_i,\l_{m-2}}}=\left\{\begin{array}{ll}
{-v^{-2}},& \ \ \textrm{if} \ i=m-3=0,\\
{v^{-4}-v^{-2}},& \ \ \textrm{if} \ i=m-2\geq1\ \textrm{or}\ i=m-3\geq1,\\
{v^{-4}},& \ \ \textrm{if} \ i=m-4,\\
{0},& \ \ \textrm{if} \ i\leq m-5\ \textrm{or}\ i=m-2=0; \end{array}
\right.$

\noindent (3) ${a_{\gamma_i,\l_{m-1}}}=\left\{\begin{array}{ll}
{-v^{-2}},& \ \ \textrm{if} \ i=m-1\ \textrm{or}\ i=m-2=0,\\
{v^{-4}-v^{-2}},& \ \ \textrm{if} \ i=m-2\geq1,\\
{v^{-4}},& \ \ \textrm{if} \ i=m-3,\\
{0},& \ \ \textrm{if} \ i\leq m-4.\end{array} \right.$

Since $\gamma_0=(m-1)x$ and $\gamma_{m-1}=2(m-1)y$ are in $X$ and
$\gamma_i=(m-1-i)x+2iy$ is in $Y$ for $i=1,2,...,m-2$, we see

\noindent (4)
${(-1)^{l(m_{\gamma_i})-l(M_{\gamma_i})}}=\left\{\begin{array}{ll}
{-1},& \ \ \textrm{if} \ i=0\ \textrm{or}\ m-1,\\
{1},& \ \ \textrm{if} \ 1\le i\le m-2\end{array}\right.$ and

\noindent (5) ${\pi_{\gamma_i}}=\left\{\begin{array}{ll}
{v+v^{-1}},& \ \ \textrm{if} \ i=0\ \textrm{or}\ m-1,\\
{1},& \ \ \textrm{if} \ 1\le i\le m-2.\end{array}\right.$

By (a), we have $b_{\lambda_{i},\l_{m-1}}=0$ if $0\leq i\leq m-3$,
$b_{\lambda_{m-1},\l_{m-1}}=1$, $b_{\lambda_{m-2},\l_{m-1}}=v^{-1}$
if $m=2$ and $b_{\lambda_{m-2},\l_{m-1}}=-v^{-2}$ if $m\geq3$. By
the proof of (a), we have
$(-1)^{l(m_{\l_i})-l(M_{\l_i})}=\pi_{\l_i}=1$ if $i\ge 1$,
$(-1)^{l(m_{\l_0})-l(M_{\l_0})}=-1$ and $\pi_{\l_0}=v+v^{-1}$.

Let $\xi=v+v^{-1}$ and $\eta=v^{-4}-v^{-2}$, $\bar\eta=v^4-v^2$. Now
using (1)--(5), Lemma 4.1 for $\gamma_i \ (0\leq i\le m-1),\
\l_{m-1}$, and 4 (a--b) we get

\noindent (6) If $m=2$, then

$-\xi \bar{b}_{\gamma_{m-1},\l_{m-1}}+v^{-1}\xi-v^{-2}=-\xi
b_{\gamma_{m-1},\l_{m-1}}+v\xi-v^{2}$,

 $-\xi
\bar{b}_{\gamma_{m-2},\l_{m-1}}+v^{-2}\xi\bar{b}_{\gamma_{m-1},\l_{m-1}}-v^{-2}=-\xi
b_{\gamma_{m-2},\l_{m-1}}+v^2\xi b_{\gamma_{m-1},\l_{m-1}}-v^{2}$;

\noindent (7) If $m\geq3$, then

 $-\xi
\bar{b}_{\gamma_{m-1},\l_{m-1}}+1-v^{-2}=-\xi
b_{\gamma_{m-1},\l_{m-1}}+1-v^{2}$,

$
\bar{b}_{\gamma_{m-2},\l_{m-1}}+v^{-2}\xi\bar{b}_{\gamma_{m-1},\l_{m-1}}
-v^{2}\eta+\eta=b_{\gamma_{m-2},\l_{m-1}}+v^{2}\xi
b_{\gamma_{m-1},\l_{m-1}} -v^{-2}\bar\eta+\bar\eta$,

$\bar{b}_{\gamma_{m-3},\l_{m-1}}-v^{-2}\bar{b}_{\gamma_{m-2},\l_{m-1}}
-v^{2}\eta+v^{-4}=b_{\gamma_{m-3},\l_{m-1}}-v^{2}
b_{\gamma_{m-2},\l_{m-1}} -v^{-2}\bar\eta+v^{4}$ \ for \  $m\geq4$,

$\bar{b}_{\gamma_{m-4},\l_{m-1}}-v^{-2}\bar{b}_{\gamma_{m-3},\l_{m-1}}
-v^{-2}=b_{\gamma_{m-4},\l_{m-1}}-v^{2} b_{\gamma_{m-3},\l_{m-1}}
-v^{2}$ \ for \  $m\geq5$,

$\bar{b}_{\gamma_{i},\l_{m-1}}-v^{-2}\bar{b}_{\gamma_{i+1},\l_{m-1}}
=b_{\gamma_{i},\l_{m-1}}-v^{2}b_{\gamma_{i+1},\l_{m-1}} $ \ for \
$1\leq i\leq m-5$,

$-\xi\bar{b}_{\gamma_{0},\l_{m-1}}-v^{-2}\bar{b}_{\gamma_{1},\l_{m-1}}+1+v^{-4}
=-\xi b_{\gamma_{0},\l_{m-1}}-v^{2} b_{\gamma_{1},\l_{m-1}}+1+v^{4}
$  \ for \   $m=3$,

$-\xi\bar{b}_{\gamma_{0},\l_{m-1}}-v^{-2}\bar{b}_{\gamma_{1},\l_{m-1}}-v^{-2}
=-\xi b_{\gamma_{0},\l_{m-1}}-v^{2} b_{\gamma_{1},\l_{m-1}}-v^{2}$
  \ for \  $m=4$,

$-\xi\bar{b}_{\gamma_{0},\l_{m-1}}-v^{-2}\bar{b}_{\gamma_{1},\l_{m-1}}
=-\xi b_{\gamma_{0},\l_{m-1}}-v^{2} b_{\gamma_{1},\l_{m-1}}$  \ for
\ $m\geq5$.

Using 4 (b) and (6), we see $b_{\gamma_{m-1},\l_{m-1}}=0$ and
$b_{\gamma_{m-2},\l_{m-1}}=v^{-1}$ if $m=2$.

Assuming $m\geq3$. Using 4 (b) and (7), we get
$b_{\gamma_{m-1},\l_{m-1}}=v^{-1};\
b_{\gamma_{m-2},\l_{m-1}}=v^{-4}-v^{-2};\ b_{\gamma_{m-3},\l_{m-1}}$
is $-v^{-3}$ if $m=3$, is $v^{-4}$ if $m\geq4$; and
$b_{\gamma_{i},\l_{m-1}}=0$ for $0\leq i\leq m-4$.

Thus (b) is proved. \hfill$\Box$

\noindent (c)  Let $m\ge 3$ be a positive integer. Set
$\nu_i=(m-2)\a+(m-2+i)\b$ for $i=0,1,...,m-2$ and
$\l_{m-1}=m\a+(2m-1)\b$.  Then we have
$${b_{\nu_i,\lambda_{m-1}}}=\left\{\begin{array}{ll}
{-v^{-3}},& \ \ \textrm{if} \ i=m-2,\\
{v^{-5}},& \ \ \textrm{if} \ i=m-3=0,\\
{-v^{-6}},& \ \ \textrm{if} \ i=m-3\geq1,\\
{0},& \ \ \textrm{if }0\le i\le m-4.\end{array} \right.$$

\noindent{\bf Proof.} The proof is similar to the proof of
(b).\hfill$\Box$

\noindent (d)  Let $\l_{m-1}=m\a+(2m-1)\b$ and
$\lambda=n\alpha+n'\beta$, where $m,n,n'$ are positive integers
satisfying $n\le n'\le 2n$ and $1\leq n\leq m-3$.  Then we have
$b_{\lambda,\lambda_{m-1}}=0$.

\noindent{\bf Proof.} We prove it by descending induction on the
partial order $\le $ in $\Lambda_r^+$.

First we  show that $b_{\lambda,\lambda_{m-1}}=0$ for
$\lambda=(m-3)\alpha+2(m-3)\b$.

Note that $m\ge 4$. By (a)--(c), we have the following observation.

(1) $b_{\nu_{m-3},\l_{m-1}}=-v^{-6},$
$b_{\nu_{m-2},\l_{m-1}}=-v^{-3},$
 $b_{\g_{m-3},\l_{m-1}}=v^{-4},$ $b_{\g_{m-2},\l_{m-1}}=v^{-4}-v^{-2},$
$b_{\g_{m-1},\l_{m-1}}=v^{-1},$ $b_{\l_{m-2},\l_{m-1}}=-v^{-2}$,
here $\nu_i,\ \g_i,\ \l_i$ are as in (c)-(a). Moreover,
$b_{\lambda',\l_{m-1}}\neq0$ for other $\l'$ in $\L_r^+$ with
$\lambda<\lambda'\leq\l_{m-1}$.

 By Proposition 4.3, we have

 (2)
$a_{\lambda,\nu_{m-3}}=-v^{-2},$
$a_{\lambda,\nu_{m-2}}=v^{-4}-v^{-2},$ $
a_{\lambda,\g_{m-3}}=v^{-4}-v^{-6},$ $ a_{\lambda,\g_{m-2}}=v^{-4},$
$a_{\lambda,\g_{m-1}}=-v^{-6},$ $ a_{\lambda,\l_{m-2}}=v^{-8}, $ $
a_{\lambda,\lambda_{m-1}}=0$.

For $\l'=\ \g_i,\  \l_i$, the values of
$(-1)^{l(m_{\l'})-l(M_{\l'}})$ and $\pi_{\l'}$ are determined in the
proofs of (a) and (b). For $\l'=\nu_i,\ \l$, it is easy to determine
the values of $(-1)^{l(m_{\l'})-l(M_{\l'}})$ and $\pi_{\l'}$ by
using Proposition 4.2. Using (1--2), Lemma 4.1 and  4 (a--b), we get
that $b_{\lambda,\l_{m-1}}=0$ for $\lambda=(m-3)\a+2(m-3)\b$ if
$m\geq4$.

Let $0\ne \l\in\Lambda_r^+$ be such that $\l<(m-3)\a+2(m-3)\b$. We
show that $b_{\l,\l_{m-1}}=0$. The induction hypothesis says that
$b_{\lambda',\l_{m-1}}=0$ for those $\lambda<\lambda'\le
(m-3)\a+(2m-6)\b$.

If $\lambda=n\a+2n\b\in X_2$, then  $n\leq m-4,$ $m\geq5$. By
Proposition 4.3, we have

(3)
$a_{\lambda,\nu_{m-3}}=a_{\lambda,\nu_{m-2}}=a_{\lambda,\g_{m-3}}
=a_{\lambda,\g_{m-2}}=a_{\lambda,\g_{m-1}}=a_{\lambda,\l_{m-2}}=a_{\lambda,\l_{m-1}}=
0$ if $m-n\geq5$; and $a_{\lambda,\nu_{m-3}}=v^{-4},$ $
a_{\lambda,\nu_{m-2}}=-v^{-6},$ $ a_{\lambda,\g_{m-3}}=v^{-8},$ $
a_{\lambda,\g_{m-2}}=a_{\lambda,\g_{m-1}}=a_{\lambda,\l_{m-2}}=a_{\lambda,\l_{m-1}}=
0$ if $m-n=4$.

Using (1), induction hypothesis, (3), Lemma 4.1 and Proposition 4.2,
we can get that $b_{\lambda,\l_{m-1}}=0$.

Similarly, for $\l=n\a+n\b$ $(1\le n\le m-3)$ or $\l=n\a+(n+i)\b\
(1\le i\le n-1,\ 1\le n\le m-3$) we can prove that
$b_{\l,\l_{m-1}}=0$.

The proof is finished. \hfill$\Box$

A reformulation of (a)--(d) is the following result.

\noindent {\bf Proposition 5.1} Let
$\lambda''=m\alpha+(2m-1)\beta\in Y_1$ for some integer $m\geq2$.
For any $0\neq\lambda\in\Lambda_r^{+}$ such that $\lambda\leq\lambda''$, we have\\
(i) If $\lambda=n\alpha+2n\beta\in X_2$ for some $n\geq2$, then
$${b_{\lambda,\lambda''}}=\left\{\begin{array}{ll}
{v^{-1}},& \ \ \textrm{if} \ n=m-1,\\
{-v^{-3}},& \ \ \textrm{if} \ n=m-2,\\
{0},& \ \ \textrm{otherwise}.\end{array} \right.$$ (ii) If
$\lambda=n\alpha+n\beta\in X_1$, then
$${b_{\lambda,\lambda''}}=\left\{\begin{array}{ll}
{v^{-1}},& \ \ \textrm{if} \ n=1\ \textrm{or}\ 2\ \textrm{and}\ m=2,\\
{v^{-5}},& \ \ \textrm{if} \ n=1\ \textrm{and}\ m=3,\\
{-v^{-3}},& \ \ \textrm{if} \ n=2\ \textrm{and}\ m=3,\\
 {0},& \ \ \textrm{otherwise}.\end{array} \right.$$ (iii) If
$\lambda=n\alpha+(2n-1)\beta\in Y_1$, then
$${b_{\lambda,\lambda''}}=\left\{\begin{array}{ll}
{1},& \ \ \textrm{if} \ n=m,\\
{v^{-4}-v^{-2}},& \ \ \textrm{if} \ n=m-1,\\
{-v^{-6}},& \ \ \textrm{if} \ n=m-2,\\
{0},& \ \ \textrm{otherwise}.\end{array} \right.$$ (iv) If
$\lambda=x^{n}y^{2n'}=(n'+n)\alpha+(2n'+n)\beta\in Y_2$, then
$${b_{\lambda,\lambda''}}=\left\{\begin{array}{ll}
{-v^{-2}},& \ \ \textrm{if} \ n=2\ \textrm{and}\ n'=m-2,\\
{v^{-4}},& \ \ \textrm{if} \ n=2\ \textrm{and}\ n'=m-3,\\
{0},& \ \ \textrm{otherwise}.\end{array} \right.$$

\def\l{\lambda}

Similarly to (a)--(d), we have the following results (e)--(h).

\noindent (e)  Let $m$ be a positive integer and $\l_i=m\a+(m+i)\b$
for $i=0,1,...,m$. Then we have
$${b_{\lambda_i,\lambda_{m}}}=\left\{\begin{array}{ll}
{1},& \ \ \textrm{if} \ i=m,\\
{v^{-1}+v^{-3}},& \ \ \textrm{if} \ i=m-1,\\
{0},& \ \ \textrm{if} \ 0\leq i\leq m-2.\end{array} \right.$$

\noindent (f)  Let $m\ge 2$ be a positive integer. Set
$\gamma_i=(m-1)\a+(m-1+i)\b$ for $i=0,1,...,m-1$ and
$\lambda_{m}=m\a+2m\b$. Then we have
$${b_{\gamma_i,\lambda_{m}}}=\left\{\begin{array}{ll}
{v^{-4}-v^{-2}},& \ \ \textrm{if} \ i=m-1,\\
{v^{-4}+v^{-2}},& \ \ \textrm{if} \ i=m-2=0,\\
{-v^{-3}-v^{-5}},& \ \ \textrm{if} \ i=m-2\geq1,\\
{0},& \ \ \textrm{if}\ 0\leq i\leq m-3.\end{array} \right.$$

\noindent (g)  Let $m\ge 3$ be a positive integer. Set
$\nu_i=(m-2)\a+(m-2+i)\b$ for $i=0,1,...,m-2$ and $\l_{m}=m\a+2m\b$.
Then we have
$${b_{\nu_i,\l_{m}}}=\left\{\begin{array}{ll}
{-v^{-6}},& \ \ \textrm{if} \ i=m-2,\\
{0},& \ \ \textrm{if}\ 0\leq i\leq m-3.\end{array} \right.$$

\noindent (h)  Let $\lambda=m'\alpha+m''\beta$ and $\l_m=m\a+2m\b\in
X_2$, where $m,m',m''$ be positive integers such that $m'<m''<2m'$
and $1\leq m'\leq m-3$. Then we have $b_{\lambda,\lambda_m}=0$.

We can  reformulate (e)--(h) as follows:

\noindent {\bf Proposition 5.2} Let $\lambda''=m\alpha+2m\beta\in
X_2$ for some integer $m\geq2$.
For any $0\neq\lambda\in\Lambda_r^{+}$ such that $\lambda\leq\lambda''$, we have that\\
(i) If $\lambda=n\alpha+2n\beta\in X_2$ for $n\geq2$, then
$${b_{\lambda,\lambda''}}=\left\{\begin{array}{ll}
{1},& \ \ \textrm{if} \ n=m,\\
{v^{-4}-v^{-2}},& \ \ \textrm{if} \ n=m-1,\\
{-v^{-6}},& \ \ \textrm{if} \ n=m-2,\\
{0},& \ \ \textrm{otherwise}.\end{array} \right.$$\\
(ii) If $\lambda=n\alpha+n\beta\in X_1$, then
$${b_{\lambda,\lambda''}}=\left\{\begin{array}{ll}
{v^{-4}+v^{-2}},& \ \ \textrm{if} \ m=2,n=1\\
{0},& \ \ \textrm{otherwise}.\end{array} \right.$$\\
(iii) If $\lambda=n\alpha+(2n-1)\beta\in Y_1$ for $2\leq n\leq m$,
then
$${b_{\lambda,\lambda''}}=\left\{\begin{array}{ll}
{v^{-1}+v^{-3}},& \ \ \textrm{if} \ n=m,\\
{-v^{-3}-v^{-5}},& \ \ \textrm{if} \ n=m-1,\\
{0},& \ \ \textrm{otherwise}.\end{array} \right.$$\\
(iv) If $\lambda=x^{n}y^{2n'}=(n'+n)\alpha+(2n'+n)\beta\in Y_2$,
then ${b_{\lambda,\lambda''}}=0.$

\noindent {\bf Remark 5.3} In [L3, (12)], Lusztig gave a conjectural
formula on $b_{\lambda,\lambda'}$ which says it is likely that, in
the case where $\check{\alpha}_{s'}(\lambda')\geq1$ for all $s'\in
I$ and $\lambda\in\Lambda_r^+$, we have
$$b_{\lambda,\lambda'}=(-1)^{l(m_{\lambda})-l(m_{\lambda'})}\frac{1}{\pi_{\lambda}}
\sum_{w\in W_{I}}(-1)^{l(w)}\Phi(w(\lambda'-\rho)-(\lambda-\rho)).$$
By Proposition 5.1 (i), we see that this conjectural formula is not
true. By [L3, (12)], we get that
$b_{(m-1)\alpha+(2m-2)\beta,m\alpha+(2m-1)\beta}=0$, which
contradicts to Proposition 5.1 (i).

Now we can compute $\mu(u,w)$ for $u,w\in c_2$. In [L3], Lusztig
computed $\mu(rsr(tsr)^n,rsr(tsr)^m)$ for $n$ is odd and $m$ is
even.

 \noindent {\bf Proposition 5.4} Let $m,n\in\mathbb{N}$. If
$m<n$ then we have
$${\mu(rsr(tsr)^m,rsr(tsr)^n)}=\left\{\begin{array}{ll}
{1},& \ \ \textrm{if} \ n-m=1,\\
{0},& \ \ \textrm{otherwise}.\end{array} \right.$$ \noindent {\bf
Proof.} We first assume that $n=2k-1$ and $m=2k'$ for some
$k>k'\in\mathbb{N}$. Let $\lambda''=(k+1)\alpha+(2k+1)\beta$ and
$\lambda=(k'+1)\alpha+2(k'+1)\beta$. Then we get
$m_{\lambda''}=rsr(tsr)^{2k-1}=rsr(tsr)^{n}$ and
$m_{\lambda}=rsr(tsr)^{2k'}=rsr(tsr)^{m}$ by Proposition 4.2 (iv).
By Proposition 5.1 (i), we get that
$${b_{\lambda,\lambda''}}=\left\{\begin{array}{ll}
{v^{-1}},& \ \ \textrm{if} \ k'=k-1,\\
{-v^{-3}},& \ \ \textrm{if} \ k'=k-2,\\
{0},& \ \ \textrm{otherwise}.\end{array} \right.$$ Therefore, by the
fact that $\mu(m_{\lambda},m_{\lambda''})$ equals to the coefficient
of $v^{-1}$ in $b_{\lambda,\lambda''}$, we get that
$\mu(m_{\lambda},m_{\lambda''})=\mu(rsr(tsr)^m,rsr(tsr)^n)=1$ if and
only if $n-m=1$; otherwise $\mu(rsr(tsr)^m,rsr(tsr)^n)=0$.

When $n$ is even and $m$ is odd, the result can be proved similarly
by Proposition 5.2 (iii). If $m-n$ is even, then
$\mu(rsr(tsr)^m,rsr(tsr)^n)=0$ is trivial. We complete the proof of
the proposition. \hfill$\Box$\


 Let $u$ and $w$ be elements in $c_2$. By 1.1 (f), we know that  if
$\mu(u,w)\neq0$ then $u\sim_L w$ or $u\sim_R w.$ First we compute
$\mu(u,w)$ for those $u,w$ such that $u\sim_L w$ and $u\sim_R w.$
After that we deal with the case $u\sim_L w$ but $u\not\sim_R w$ and
the case $u\not\sim_L w$ but $u\sim_R w.$

\noindent {\bf Lemma 5.5} Let $m,n\in\mathbb{N}$. If $m<n$ then we
have
$${\mu(r(tsr)^m,rsr(tsr)^n)}=\left\{\begin{array}{ll}
{1},& \ \ \textrm{if} \ n-m=1,\\
{0},& \ \ \textrm{otherwise}.\end{array} \right.$$ \noindent {\bf
Proof.} We have two left strings with respect to $S'=\{r,s\}:$
$$u_1=r(tsr)^m,\ u_2=sr(tsr)^m,\ u_3=rsr(tsr)^m$$
and$$w_1=r(tsr)^n,\ w_2=sr(tsr)^n,\ w_3=rsr(tsr)^n.$$ By 1.1 (g), we
get that $$\mu(r(tsr)^m,\ rsr(tsr)^n)=\mu(rsr(tsr)^m,\ r(tsr)^n).$$
Since $t\in \mathcal {L}(r(tsr)^n)\backslash \mathcal
{L}(rsr(tsr)^m)$, then we get the result using 1.1 (b).\hfill$\Box$\

\noindent {\bf Lemma 5.6} Let $m,n\in\mathbb{N}$. If $m<n$ then
 we have
$${\mu(r(tsr)^mts,\ rsr(tsr)^nts)}=\left\{\begin{array}{ll}
{1},& \ \ \textrm{if} \ n-m=1,\\
{0},& \ \ \textrm{otherwise}.\end{array} \right.$$

\noindent {\bf Proof.} The proof is very similar to that of Lemma
5.5\hfill$\Box$\

For any two integers $m,n,$ we get two left strings with respect to
$S'=\{r,s\}:$
$$u_1=r(tsr)^m,\ u_2=sr(tsr)^m,\ u_3=rsr(tsr)^m$$
and$$w_1=r(tsr)^n,\ w_2=sr(tsr)^n,\ w_3=rsr(tsr)^n.$$ By 1.1 (g), we
have
$$\mu(r(tsr)^m,r(tsr)^n)=\mu(rsr(tsr)^m,rsr(tsr)^n)$$
and
$$\mu(sr(tsr)^m,sr(tsr)^n)=\mu(r(tsr)^m,r(tsr)^n)+\mu(r(tsr)^m,rsr(tsr)^n).$$
Similarly, we get
$\mu(sr(tsr)^mts,sr(tsr)^nts)$\begin{eqnarray*}&=&\mu(r(tsr)^mts,r(tsr)^nts)+\mu(r(tsr)^mts,rsr(tsr)^nts)\\
&=&\mu(rsr(tsr)^mts,rsr(tsr)^nts)+\mu(r(tsr)^mts,rsr(tsr)^nts)\\
&=&\mu((rst)^{m+1}rs,(rst)^{n+1}rs)+\mu(r(tsr)^mts,rsr(tsr)^nts).
\end{eqnarray*}
By 1.1 (a), we have
$$\mu((rst)^{m+1}rs,(rst)^{n+1}rs)=\mu(sr(tsr)^{m+1},sr(tsr)^{n+1}).$$
We also get the following identities:
\\$\mu(r(tsr)^mt,r(tsr)^nt)$
\begin{eqnarray*}&=&\mu(rsr(tsr)^mt,rsr(tsr)^nt)=\mu((rst)^{m+1}r,(rst)^{n+1}r)\\
&=&\mu(r(tsr)^{m+1},r(tsr)^{n+1})=\mu(rsr(tsr)^{m+1},rsr(tsr)^{n+1}),
\end{eqnarray*}
\begin{eqnarray*}\mu(rsr(tsr)^mtst,rsr(tsr)^ntst)&=&\mu(rsr(tsr)^mt,rsr(tsr)^nt)\\&=&\mu(r(str)^{m+1},r(str)^{n+1})\\
&=&
\mu((rst)^{m+1}r,(rst)^{n+1}r)\\&=&\mu(r(tsr)^{m+1},r(tsr)^{n+1})\end{eqnarray*}(regarded
as right strings with respect to $S'=\{s,t\}$),
\begin{eqnarray*}\mu((rts)^{m+1}t,(rts)^{n+1}t)
&=&\mu(r(tsr)^mtst,r(tsr)^ntst\\
&=&\mu(r(tsr)^mt,r(tsr)^nt),
\end{eqnarray*}
\begin{eqnarray*}\mu(sr(tsr)^mt,sr(tsr)^nt)&=&\mu(st(rst)^mr,st(rst)^nr)\\
&=&\mu(st(rst)^mrsr,st(rst)^nrsr)\\
&=&\mu(sr(tsr)^{m+1},sr(tsr)^{n+1}),
\end{eqnarray*}
$$\mu(sr(tsr)^mtst,sr(tsr)^ntst)=\mu(sr(tsr)^mt,sr(tsr)^nt)$$and
\begin{eqnarray*}\mu((tsr)^mtst,(tsr)^ntst)&=&\mu((tsr)^mt,(tsr)^nt)\\
&=&\mu(t(str)^m,t(str)^n)\\
&=&\mu((rts)^mt,(rts)^nt).
\end{eqnarray*}

With these identities, Proposition 5.4, Lemma 5.5 and Lemma 5.6, we
obtain the following theorem.

\noindent {\bf Theorem 5.7} For elements $u, w\in c_2$ such that
$u<w$, $u\sim_L w$ and
$u\sim_R w,$ we have \\
(i) If $s$ is not in $\mathcal {L}(u)\cup \mathcal {R}(u)$, then
$${\mu(u,w)}=\left\{\begin{array}{ll}
{1},& \ \ \textrm{if} \ l(w)-l(u)=3,\\
{0},& \ \ \textrm{otherwise}.\end{array} \right.$$ \\
(ii) If $s\in \mathcal {L}(u)\setminus \mathcal {R}(u)$ or $s\in
\mathcal {R}(u)\setminus \mathcal {L}(u),$ then
$${\mu(u,w)}=\left\{\begin{array}{ll}
{2},& \ \ \textrm{if} \ l(w)-l(u)=3,\\
{0},& \ \ \textrm{otherwise}.\end{array} \right.$$ \\
(iii) If $s\in \mathcal {L}(u)\cap \mathcal {R}(u),$ then
$${\mu(u,w)}=\left\{\begin{array}{ll}
{3},& \ \ \textrm{if} \ l(w)-l(u)=3,\\
{0},& \ \ \textrm{otherwise}.\end{array} \right.$$

\noindent {\bf Proof.} The only thing to notice is the fact: for any
$u,w\in c_2,$ $u\sim_L w\Leftrightarrow \mathcal {R}(u)=\mathcal
{R}(w)$ and $u\sim_R w\Leftrightarrow \mathcal {L}(u)=\mathcal
{L}(w).$\hfill$\Box$\

 In the following we compute $\mu(u,w)$ for those  $u, w\in c_2$  satisfying
 (1) $u\leq w$, (2) $u\sim_L w$ and $u\nsim_R w,$ or $u\sim_R w$
and $u\nsim_L w.$

\noindent {\bf Theorem 5.8} Let $u,w\in c_2$ with $u\leq w$. Assume
that
 $u\sim_Lw$ and
$u\nsim_Rw.$ We have \\
(i) If $\mathcal {L}(w)\nsubseteq \mathcal {L}(u),$ then
$${\mu(u,w)}=\left\{\begin{array}{ll}
{1},& \ \ \textrm{if} \ l(w)-l(u)=1,\\
{0},& \ \ \textrm{otherwise}.\end{array} \right.$$\\
(ii) If $\mathcal {L}(w)\subseteq \mathcal {L}(u),$ then
$${\mu(u,w)}=\left\{\begin{array}{ll}
{1},& \ \ \textrm{if} \ l(w)-l(u)=5,\\
{0},& \ \ \textrm{otherwise}.\end{array} \right.$$

\noindent {\bf Proof.} Obviously, we have the fact: for any $u,w\in
c_2,$ $u\sim_L w\Leftrightarrow \mathcal {R}(u)=\mathcal {R}(w)$ and
$u\sim_R w\Leftrightarrow \mathcal {L}(u)=\mathcal {L}(w).$ By the
assumption on $u,w,$ we get $R(u)=R(w)\ \textrm{but}\ L(u)\neq
L(w).$ The proof of (i) is trivial by 1.1(b).

Now we assume that $L(w)\subseteq L(u),$ thus $L(u)=\{r,t\},\
L(w)=\{r\}$ or $\{t\}.$ Without loss of generality, we can assume
that $L(w)=\{r\}.$ Then $$u\in \{rt(str)^m,\ rt(str)^ms,\
rt(str)^msr,\ rt(str)^mst\mid m\in \mathbb{N}\}$$
$$w\in \{r(str)^n,\ r(str)^ns,\
r(str)^nsr,\ r(str)^nst\mid n\in \mathbb{N}\}$$ If
$u=rt(str)^m=r(tsr)^mt,$ then $w=r(str)^n=rsr(tsr)^{n-1}t.$ By 1.1
(g), we get
$\mu(r(tsr)^mt,rsr(tsr)^{n-1}t)=\mu(rsr(tsr)^mt,r(tsr)^{n-1}t).$
Since $t\in L(r(tsr)^{n-1}t)\backslash L(rsr(tsr)^mt)$, then by 1.1
(b) we get that
$${\mu(rsr(tsr)^mt,r(tsr)^{n-1}t)}=\left\{\begin{array}{ll}
{1},& \ \ \textrm{if} \ n-m=2,\\
{0},& \ \ \textrm{otherwise}.\end{array} \right.$$\\
Thus we get
$${\mu(rt(str)^m,r(str)^n)}=\left\{\begin{array}{ll}
{1},& \ \ \textrm{if} \ n-m=2,\\
{0},& \ \ \textrm{otherwise}.\end{array} \right.$$ Similarly, we
have
$${\mu(rt(str)^ms,r(str)^ns)}=\left\{\begin{array}{ll}
{1},& \ \ \textrm{if} \ n-m=2,\\
{0},& \ \ \textrm{otherwise};\end{array} \right.$$
$${\mu(rt(str)^msr,r(str)^nsr)}=\left\{\begin{array}{ll}
{1},& \ \ \textrm{if} \ n-m=2,\\
{0},& \ \ \textrm{otherwise}\end{array} \right.$$and
$${\mu(rt(str)^mst,r(str)^nst)}=\left\{\begin{array}{ll}
{1},& \ \ \textrm{if} \ n-m=2,\\
{0},& \ \ \textrm{otherwise}.\end{array} \right.$$ By the
assumption, we know that in this case $n-m=2$ is equivalent to
$l(w)-l(u)=5$, thus (ii) holds. \hfill$\Box$\

Similarly, we have the following theorem.

\noindent {\bf Theorem 5.9} Let $u,w\in c_2$ with $u\leq w$. Assume
that $u\sim_Rw$ and
$u\nsim_Lw.$ We have \\
(i) If $R(w)\nsubseteq R(u),$ then
$${\mu(u,w)}=\left\{\begin{array}{ll}
{1},& \ \ \textrm{if} \ l(w)-l(u)=1,\\
{0},& \ \ \textrm{otherwise}.\end{array} \right.$$\\
(ii) If $R(w)\subseteq R(u),$ then
$${\mu(u,w)}=\left\{\begin{array}{ll}
{1},& \ \ \textrm{if} \ l(w)-l(u)=5,\\
{0},& \ \ \textrm{otherwise}.\end{array} \right.$$


\section {The proof of Theorem 2.2}
Now we can give the proof of Theorem 2.2. The only thing we need to
prove is: for some $d\in {\cal D}_0,$ if $z\sim_Lz^{-1}\sim_Ld$ and
$\tilde{\mu}(z,d)\neq0,$ then $z\in {\cal D}_1.$

If $d\in c_1$ and $\tilde{\mu}(z,d)\neq0,$ then by the proof of
Theorem 3.3 we know that $z\sim_Ld$ and $z^{-1}\sim_Ld$ cannot hold
at the same time. Also, if $d\in c_0$ and $\tilde{\mu}(z,d)\neq0,$
we know that $z\sim_Ld$ and $z^{-1}\sim_Ld$ cannot hold at the same
time, too (see the proof of Theorem 3.1 in [SX]).

Now we assume that $d\in c_2\ \textrm{and}\  z\sim_Ld\sim_Lz^{-1}.$
Then we get $L(z)=R(z)=L(d).$ By Theorem 5.7, those $z$ satisfying
$z\sim_Lz^{-1}\sim_Ld$ and $\tilde{\mu}(d,z)\neq0$ for some
$d\in{\cal D}_0\cap c_2$ are just $rtstr, strstrs, rstrstrsr,
tstrstrsr$. By some computations, we can check that these four
elements are all in ${\cal D}_1.$ The values of the
$\textbf{a}$-function and the length function on these elements are
obvious. We just need to get the degree $\delta (w)$ of polynomials
$P_{e,w}$ for these $w$ by using the formula in [KL1, (2.2.c)]. For
example, we have
$P_{e,rtstr}=P_{r,rtstr}=P_{tr,rtstr}=P_{r,rstr}+qP_{tr,rstr}=q+1$.
Thus we have $l(rtstr)=5,\ a(rtstr)=2$ and $\delta (rtstr)=1$.
Therefore $rtstr\in {\cal D}_1$. With the same method, we can check
that $strstrs, rstrstrsr, tstrstrsr$ are all in ${\cal D}_1$.
Therefore, Conjecture 2.1 holds for an affine Weyl group of type
$\tilde{B}_2$.

\noindent {\bf Remark 6.1} In fact, we can show that the set ${\cal
D}_1$ is finite for an affine Weyl group of type $\tilde{B}_2$ by
Theorem 3.3, Theorem 5.7 and some results about $\mu(u,w)$ for $u<w$
and $(u,w)\in c_0\times c_0$. For more details see [W].

\section {Computing $\mu(u,w)$ for $a(u)< a(w)$}
In this section we will compute $\mu(u,w)$ clearly for those $u<w$
such that $a(u)<a(w)$. We will see that in this case
$\mu(u,w)\leq1$.

\noindent {\bf Theorem 7.1} For $u<w$ such that $a(u)=1\
\textrm{or}\ 2\ \textrm{and}\ a(w)=4,$ we have
$${\mu(u,w)}=\left\{\begin{array}{ll}
{1},& \ \ \textrm{if} \ l(w)-l(u)=1,\\
{0},& \ \ \textrm{otherwise}.\end{array} \right.$$

\noindent {\bf Proof.} By 1.1 (f), if $\mu(u,w)\neq0$ we have
$w\leq_L u$ and $w\leq_R u.$ Then by 1.1 (d) and (e), if
$\mu(u,w)\neq0$ we have $R(u)\subseteq R(w)$ and $L(u)\subseteq
L(w).$ If $R(u)\subsetneq R(w)$ or $L(u)\subsetneq L(w),$ the result
is obvious by 1.1(b) and (c). In the following we assume that
$R(u)=R(w)$ and $L(u)=L(w).$

(1) Assume that $a(u)=1.$ If $L(u)=L(w)=\{s\},$ we have
$L(sw)=\{r,t\},\ L(su)=\{r\}\ \textrm{or}\ \{t\}\ \textrm{or}\
\varnothing.$ Then there exists some $s'\in S$ such that $s'\in
L(sw)\backslash L(su).$ By [KL1, (2.2.c)], we have
$$P_{u,w}=P_{su,sw}+qP_{u,sw}-\sum_{{u\leq z\prec sw}\atop {sz\leq
z}}\mu(z,sw)q^{\frac{1}{2}(l(sw)-l(z)+1)}P_{u,z}.$$ Also we have
$P_{su,sw}=P_{s'su,sw}$ and $P_{u,sw}=P_{rtu,sw}.$ If $s'su=sw$, we
have $\mu(u,w)=1$. If $l(w)-l(u)>1$ we get
$\textrm{deg}P_{su,sw}<\frac{1}{2}(l(w)-l(u)-1)$ and
$\textrm{deg}(qP_{u,sw})<\frac{1}{2}(l(w)-l(u)-1).$ By the
non-negativity of the coefficients of Kazhdan--Lusztig polynomials,
we get $\mu(u,w)=0$ if $l(w)-l(u)\neq1.$

If $L(u)=L(w)=\{r\},$ we can assume that $L(ru)\neq\varnothing$,
since if $u=r$ the proof is similar to the above. We have
$L(ru)=L(rw)=\{s\}$ and $a(ru)=1,\ a(rw)=4.$ By [KL1, (2.2.c)], we
have
$$P_{u,w}=P_{ru,rw}+qP_{u,rw}-\sum_{{u\leq z\prec rw}\atop{rz\leq
z}}\mu(z,rw)q^{\frac{1}{2}(l(rw)-l(z)+1)}P_{u,z}.$$ If
$l(w)-l(u)>1,$ by the above proof, we get
$\textrm{deg}P_{ru,rw}<\frac{1}{2}(l(w)-l(u)-1).$ Moreover,
$P_{u,rw}=P_{su,rw},\ L(su)=L(rw)=\{s\}.$ But $l(rw)-l(su)\neq1,$
then we get $\textrm{deg}(qP_{u,rw})<\frac{1}{2}(l(w)-l(u)-1)$ by
the above proof. Thus $\mu(u,w)=0$ in this case.

If $L(u)=L(w)=\{t\},$ the proof is similar to the case
$L(u)=L(w)=\{r\}$.

(2) Assume that $a(u)=2.$ The proof is similar to (1). The only
thing to note is the following facts: If $L(u)=L(w)=\{r,t\},$ then
there exists some $s'\in\{r,t\}$ such that $L(s'w)=\{s,s''\},$ where
$s'\neq s''\in\{r,t\},$ and $L(s'u)=\{s''\}.$ We have
$P_{u,s'w}=P_{s''su,s'w},$ where $s''su\geq su\geq u.$ Then we
consider the identity
$$P_{u,w}=P_{s'u,s'w}+qP_{u,s'w}-\sum_{{u\leq z\prec s'w}\atop{s'z\leq
z}}\mu(z,s'w)q^{\frac{1}{2}(l(s'w)-l(z)+1)}P_{u,z}.$$

If $L(u)=L(w)=\{s\},$ we have $L(su)=L(sw)=\{r,t\}.$ This case is
trivial.

If $L(u)=L(w)=\{r\},$ we have $L(ru)=L(rw)=\{s\}.$ If $l(w)-l(u)>1,$
we get $P_{u,rw}=P_{su,rw}$ and $L(su)=\{r,s\}.$ We consider the
identity
$$P_{u,w}=P_{ru,rw}+qP_{u,rw}-\sum_{{u\leq z\prec rw}\atop{rz\leq
z}}\mu(z,rw)q^{\frac{1}{2}(l(rw)-l(z)+1)}P_{u,z}.$$

If
$\mu(su,rw)\neq0,$ we take $z=su$ in the sum, then the item\\
$\mu(su,rw)q^{\frac{1}{2}(l(w)-l(u)-1)}$ appearing in the sum is
killed by $qP_{u,rw}.$

If $L(u)=L(w)=\{t\},$ the proof is similar to the case
$L(u)=L(w)=\{r\}.$ We complete the proof . \hfill$\Box$\

We define a subset $E$ of $c_1\times c_2$ by\\
$E=\{(st,srtst),(sr,srtsr),(r,rsrtsr),
(t,tsrtst),(rst,rsrtst),(tsr,tsrtsr)\}.$ We have the following
theorem.

\noindent {\bf Theorem 7.2} For $u<w$ such that $a(u)=1\
\textrm{and}\ a(w)=2,$ we have
$${\mu(u,w)}=\left\{\begin{array}{ll}
{1},& \ \ \textrm{if} \ l(w)-l(u)=1 \ \textrm{or}\ (u,w)\in E,\\
{0},& \ \ \textrm{otherwise}.\end{array} \right.$$

\noindent {\bf Proof.} By 1.1 (f), if $\mu(u,w)\neq0$ we get
$R(u)\subseteq R(w)$ and $L(u)\subseteq L(w).$ If $R(u)\subsetneq
R(w)$ or $L(u)\subsetneq L(w)$ the result is obvious by 1.1(b) and
(c). In the following we assume that $L(u)=L(w)=\{s'\}$ and
$R(u)=R(w)=\{s''\}.$

(1) If $s'=s,s''=r.$ We have $u=(stsr)^m\ \textrm{or}\
u=sr(stsr)^{m'},w=(str)^nsr$ for some $1\leq m,n\in \mathbb{N},0\leq
m'\in \mathbb{N}.$ If $u=(stsr)^m,w=(str)^nsr$ for some $1\leq
m,n\in \mathbb{N},$ we have $L(su)=\{t\},L(sw)=\{r,t\}.$ Consider
the identity
$$P_{u,w}=P_{su,sw}+qP_{u,sw}-\sum_{{u\leq z\prec sw}\atop{sz\leq
z}}\mu(z,sw)q^{\frac{1}{2}(l(sw)-l(z)+1)}P_{u,z}.$$ We get
$P_{su,sw}=P_{rsu,sw}$ and $P_{u,sw}=P_{rtu,sw}.$ If $l(w)-l(u)=1$
(i.e. $ m=n=1$) we get $\mu(u,w)=1.$ Otherwise, we have
$\textrm{deg}P_{su,sw}<\frac{1}{2}(l(w)-l(u)-1),\
\textrm{deg}(qP_{u,sw})<\frac{1}{2}(l(w)-l(u)-1),$ thus
$\mu(u,w)=0.$

If $u=sr(stsr)^{m'},w=(str)^nsr$ for some $0\leq m'\in
\mathbb{N},1\leq n\in \mathbb{N},$ we have
$L(su)=\{r\},L(sw)=\{r,t\}.$ We get $P_{su,sw}=P_{tsu,sw}$ and
$P_{u,sw}=P_{rtu,sw}.$ In this case there is no $m',n\in \mathbb{N}$
satisfying $tsu=sw,$ thus
$\textrm{deg}P_{su,sw}<\frac{1}{2}(l(w)-l(u)-1).$ Also we get
$P_{rtu,sw}=1$ if $m'=0,n=1.$ Otherwise,
$\textrm{deg}(qP_{u,sw})<\frac{1}{2}(l(w)-l(u)-1).$ Then we get
$\mu(u,w)=1$ if $u=sr,w=srtsr;$ otherwise $\mu(u,w)=0.$

($1'$) The proofs for cases (i) $s'=s,s''=t;$ (ii) $s'=r,s''=s$ and
(iii) $s'=t,s''=s$ are similar to (1).

(2) If $s'=r,s''=t.$ We have $u=r(stsr)^mst,w=r(str)^nst$ for some
$0\leq m\in \mathbb{N}$ and $1\leq n\in \mathbb{N}.$ Then
$L(ru)=L(rw)=\{s\}.$ Consider the identity
$$P_{u,w}=P_{ru,rw}+qP_{u,rw}-\sum_{{u\leq z\prec rw}\atop{rz\leq
z}}\mu(z,rw)q^{\frac{1}{2}(l(rw)-l(z)+1)}P_{u,z}.$$ By ($1'$), we
know that $\mu(ru,rw)=1$ if $m=0$ and $n=1$; $\mu(ru,rw)=0$
otherwise. Moreover, $P_{u,rw}=P_{su,rw}$ with $su>u$ and $a(su)=1.$
By ($1'$), we know that $\mu(su,rw)=1$ if $m=0,n=1;$ otherwise
$\mu(su,rw)=0.$ Then we get that $\mu(u,w)=1$ if $u=rst,w=rstrst;$
otherwise $\mu(u,w)=0.$

($2'$) The case $s'=t,s''=r$ is similar to (2).

(3) If $s'=s''=r.$ We get $u=r(stsr)^m,w=r(str)^nsr$ for some $0\leq
m\in \mathbb{N},1\leq n\in \mathbb{N}.$ If $m=0$ (i.e $u=r$) we get
$\mu(r,w)=1$ if $w=rsrtsr;$ otherwise $\mu(r,w)=0.$ In the following
we assume that $m\geq1,$ thus we get $L(ru)=L(rw)=\{s\}.$ By (1), we
know that $\mu(ru,rw)=1$ if $m=n=1;$ otherwise $\mu(ru,rw)=0.$
Moreover, we have $P_{u,rw}=P_{su,rw}$ with $su>u$ and $a(su)=1.$ By
($1'$), we know that $\mu(su,rw)=1$ if $m=0,n=1;$ otherwise
$\mu(su,rw)=0.$ Thus we get the result.

($3'$) The case $s'=s''=t$ is similar to (3).

(4) If $s'=s''=s$ we shall prove that $\mu(u,w)=1$ if $l(w)-l(u)=1,$
otherwise, $\mu(u,w)=0.$

We have
$u\in\{(stsr)^ms,(stsr)^msts,sr(stsr)^ms,sr(stsr)^msts\mid0\leq m\in
\mathbb{N}\}$ and $w=(str)^ns$ for some $1\leq n\in \mathbb{N}.$ If
$u=s$ we can compute that $\mu(s,w)=0$ easily. In the following we
assume that $u\neq s.$

If $u=(stsr)^ms,w=(str)^ns$ for some $1\leq m,n\in \mathbb{N},$ we
have $L(su)=\{t\},L(sw)=\{r,t\}.$ Then we get $\mu(su,sw)=1$ if
$l(w)-l(u)=1,$ otherwise, $\mu(su,sw)=0.$ Also, we have
$P_{u,sw}=P_{rtu,sw}\neq1,$ thus
$\textrm{deg}(qP_{u,sw})<\frac{1}{2}(l(w)-l(u)-1).$ So we get
$\mu(u,w)=1$ if $l(w)-l(u)=1,$ otherwise, $\mu(u,w)=0.$

If $u\in\{(stsr)^msts,sr(stsr)^ms,sr(stsr)^msts\mid0\leq
m\in\mathbb{N}\},$ the proof is similar. We complete the proof of
the theorem. \hfill$\Box$\

\section { Computing $\mu(u,w)$ for $(u,w)\in c_2\times c_1$}
In this section we compute $\mu(u,w)$ clearly for $a(u)=2\
\textrm{and}\ a(w)=1.$ We'll see that $\mu(u,w)\leq1$ in this case.

\noindent {\bf Proposition 8.1} For any $1\leq n,m\in \mathbb{N},$
we have
$${\mu(rsr(tsr)^n,(rsts)^mr)}=\left\{\begin{array}{ll}
{1},& \ \ \textrm{if} \ n=1,m\geq2 ,\\
{0},& \ \ \textrm{otherwise}.\end{array} \right.$$

\noindent {\bf Proof.} The proof is similar to the proof of
Proposition 5.4. In the following we assume that $n$ is odd,
otherwise $\mu(rsr(tsr)^n,(rsts)^mr)=0$. Let
$\lambda=\frac{n+3}{2}\alpha+(n+2)\beta=xy^{n+1}\ \textrm{and}\
\lambda''=(m+1)\alpha+(m+1)\beta=x^{m+1},$ then we have
$m_\lambda=rsr(tsr)^n\ \textrm{and}\ m_{\lambda''}=(rsts)^mr,$ where
$m_\lambda\ \textrm{and}\ m_{\lambda''}$ are the unique elements of
minimal length in $W_0\lambda W_0$ and $W_0\lambda'' W_0$
respectively. By Section 4.1, we know that
$\mu(rsr(tsr)^n,(rsts)^mr)=\textrm{Res}_{v=0}b_{\lambda,\lambda''}.$
Due to Lemma 4.1, Proposition 4.2 and Proposition 4.3, we can get
that
$${b_{\lambda,\lambda''}}=\left\{\begin{array}{ll}
{v^{-1}+v^{-3}},& \ \ \textrm{if} \ n=1,m=2 ,\\
{v^{-1}-v^{-5}},& \ \ \textrm{if} \ n=1,m\geq3,\\
{0},& \ \ \textrm{otherwise}.\end{array} \right.$$\\
Then the result is obvious. \hfill$\Box$\

Following this proposition we can get our main result in this
section.

\noindent {\bf Theorem 8.2} For $u<w$ satisfying $a(u)=2$ and $a(w)=1,$ we have\\
(i) If $L(w)=R(w)=\{s\}$, then
$${\mu(u,w)}=\left\{\begin{array}{ll}
{1},& \ \ \textrm{if} \ u=srts ,\\
{0},& \ \ \textrm{otherwise}.\end{array} \right.$$\\
(ii) If $L(w)=\{s\},R(w)=\{s'\},$ where $s'=r$ or $t,$ then
$${\mu(u,w)}=\left\{\begin{array}{ll}
{1},& \ \ \textrm{if} \ u=srtss'\ \textrm{or}\ srt ,\\
{0},& \ \ \textrm{otherwise}.\end{array} \right.$$\\
 If $L(w)=\{s'\},R(w)=\{s\},$ where $s'=r$ or $t,$ then
$${\mu(u,w)}=\left\{\begin{array}{ll}
{1},& \ \ \textrm{if} \ u=s'srts\ \textrm{or}\ rts ,\\
{0},& \ \ \textrm{otherwise}.\end{array} \right.$$\\
(iii) If $L(w)=R(w)=\{s'\},$ where $s'=r$ or $t,$ $w\neq tsrst$ and
$rstsr,$ then
$${\mu(u,w)}=\left\{\begin{array}{ll}
{1},& \ \ \textrm{if} \ u\in \{s'srtss',rtss',s'srt,rt\} ,\\
{0},& \ \ \textrm{otherwise}.\end{array} \right.$$\\
If $L(w)=R(w)=\{s'\},$ where $s'=r$ or $t,$ $w= tsrst$ or $rstsr,$
then
$${\mu(u,w)}=\left\{\begin{array}{ll}
{1},& \ \ \textrm{if} \ u\in \{rtss',s'srt,rt\} ,\\
{0},& \ \ \textrm{otherwise}.\end{array} \right.$$\\
(iv) If $L(w)=\{s'\},R(w)=\{s''\},$ where $s'\neq s''\in \{r,t\}$,
$w\neq rst(srst)^2$ and $tsr(stsr)^2,$ then
$${\mu(u,w)}=\left\{\begin{array}{ll}
{1},& \ \ \textrm{if} \ u\in \{s'srtss'',rtss'',s'srt,rt\} ,\\
{0},& \ \ \textrm{otherwise}.\end{array} \right.$$\\
If $L(w)=\{s'\},R(w)=\{s''\},$ where $s'\neq s''\in \{r,t\}$, $w=
rst(srst)^2$ or $tsr(stsr)^2,$  then
$${\mu(u,w)}=\left\{\begin{array}{ll}
{1},& \ \ \textrm{if} \ u\in \{rtsrt,s'srtss'',rtss'',s'srt,rt\} ,\\
{0},& \ \ \textrm{otherwise}.\end{array} \right.$$

\noindent {\bf Proof.} For any $u<w$ satisfying $a(u)=2$ and
$a(w)=1,$ if $\mu(u,w)\neq0,$ by 1.1 (d, e, f) we get $L(w)\subseteq
L(u)\ \textrm{and}\ R(w)\subseteq R(u).$

We prove (iii) first. We assume that $L(w)=R(w)=\{r\},$ then
$w=(rsts)^mr$ for some integer $m\geq1$ and $L(u),\ R(u)=\{r\}$ or
$\{r,t\}$. If $m=1$, the claim is very easy to prove. In the
following we assume that $m\geq2.$ Consider the two left strings
w.r.t $\{r,s\}$: $u_1=r(tsr)^n,u_2=sr(tsr)^n,u_3=rsr(tsr)^n$ and
$w_1=s(tsrs)^{m-1}tsr=(stsr)^m,w_2=rs(tsrs)^{m-1}tsr=(rsts)^mr,w_3
=srs(tsrs)^{m-1}tsr=s(rsts)^mr=(srst)^msr.$

By 1.1 (g), we get $\mu(u_1,w_2)=\mu(u_3,w_2)$, i.e.
$\mu(r(tsr)^n,(rsts)^mr)=\mu(rsr(tsr)^n,(rsts)^mr).$ Similarly, we
get\\ $\mu(rt(str)^n,(rsts)^mr=\mu(r(str)^{n+1},(rsts)^mr)$ (left
strings w.r.t $\{r,s\}$). When $L(w)=R(w)=\{t\},$ we can also get
the similar result. Using the facts $\mu(u,w)=\mu(u^{-1},w^{-1})$
and Proposition 8.1, we can get (iii).

We now prove (ii). Assume that $L(w)=\{s\}$ and $R(w)=\{r\},$ then\\
$w\in\{(srst)^msr,(stsr)^m\mid m\geq1\}.$ By 1.1 (g), we get\\
$\mu((str)^nsr,(stsr)^m)=\mu((str)^nsr,(srst)^msr)=\mu(rsr(tsr)^n,(rsts)^mr)$
(left strings w.r.t $\{r,s\}$) and
$\mu((str)^nsr,(srst)^msr)=\mu((str)^n,(srst)^msr)\\=\mu((str)^nsr,(srst)^msr)$
(right strings w.r.t $\{r,s\}$). Then (ii) holds with the same
reason of (iii).

The proofs of (i) and (iv) are very similar to those of (ii) and
(iii).\hfill$\Box$\

\noindent {\bf Remark 8.3} (1) From Theorem 8.2, we see that the
W-graph of type $\tilde{B}_2$ is nonlocally finite. More
specifically, we have some element $u$ in $W$ such that there are
infinitely many elements $w$ in $W$ satisfying $\mu(u,w)\neq0$. This
is one of the main interests of [L3] (the definition of W-graph can
be found in [KL1]).

 (2) For those $u<w$ satisfying $a(u)=4$ and
$a(w)=1$ or 2, we can only get part of the leading coefficients. For
example, we have $\mu(u,w)=0$ when $u$ and $w$ are both of the
minimal length in their double cosets $W_0uW_0$ and $W_0wW_0$,
respectively. As a consequence, we determine all $\mu(u,w)$ for
those $u<w$ satisfying $l(w)-l(u)\leq3$.


\section*{Acknowledgment}
The results of this paper were obtained during my PhD. studies at
Chinese Academy of Sciences. I would like to express deep gratitude
to my supervisor Professor N.Xi, whose guidance and support were
crucial for the successful completion of this paper. Also I must
thank Professor G.Lusztig for his very helpful answers to my
questions about this paper. I am very grateful to the referee for
careful reading and valuable comments, which improve the paper
significantly.

\providecommand{\bysame}{\leavevmode\hbox
to3em{\hrulefill}\thinspace}
\providecommand{\MR}{\relax\ifhmode\unskip\space\fi MR }
\providecommand{\MRhref}[2]{%
  \href{http://www.ams.org/mathscinet-getitem?mr=#1}{#2}
} \providecommand{\href}[2]{#2}

\end{document}